\documentclass[11pt]{amsart}
\usepackage{amsmath,amssymb,amsthm,mathtools}
\usepackage[expansion=false]{microtype}
\usepackage{needspace}
\usepackage{enumitem,booktabs,array,longtable}
\usepackage[hmargin=1.05in,vmargin=0.96in]{geometry}
\usepackage[svgnames]{xcolor}
\usepackage{hyperref}
\hypersetup{colorlinks=true,allcolors=Navy,breaklinks=true,
  pdftitle={A complete classification of permutation binomials of the form X\textasciicircum r(X\textasciicircum\{2(q-1)\}+a) over finite fields},
  pdfauthor={Zhichao Tang and Xiang Fan}}
\allowdisplaybreaks[2]
\numberwithin{equation}{section}
\newtheorem{theorem}{Theorem}[section]
\newtheorem{proposition}[theorem]{Proposition}
\newtheorem{lemma}[theorem]{Lemma}
\newtheorem{corollary}[theorem]{Corollary}
\theoremstyle{remark}

\title[Permutation binomials $X^r(X^{2(q-1)}+a)$]
{A complete classification of permutation binomials of the form
$X^r(X^{2(q-1)}+a)$ over finite fields}
\author{Zhichao Tang}
\author{Xiang Fan}
\address{School of Mathematics, Sun Yat-sen University, Guangzhou 510275, China}
\email[Xiang Fan]{fanx8@mail.sysu.edu.cn}
\begin{document}
\begin{abstract}
Let $q$ be a prime power. We determine, for all positive integers $e$
and $r$, the coefficients $a\in\mathbb{F}_{q^e}$ for which
$X^r(X^{2(q-1)}+a)$ permutes $\mathbb{F}_{q^e}$. For odd $q$ and
$a\ne0$ there are no such permutations when $e\geqslant4$ is even. Over
$\mathbb{F}_{q^2}$, with $r$ unrestricted, the permutation binomials are
characterized by congruence and coprimality conditions on $r$ and
equations of degree at most two for $(-a)^{(q+1)/2}$, together with one
sporadic family over $\mathbb{F}_{11^2}$; this gives an explicit form of
Hou's finiteness theorem for each fixed odd $r>3$. The remaining cases
follow from the known classification for the exponent difference $q-1$.
In the quadratic case, power sums whose binomial expansions are sparse
by Lucas' theorem determine the exponent one $p$-adic digit at a time.
In even degree at least four, an exponent congruence reduces the problem
to the existence of a collision, which we obtain by counting norm ratios
on the projective line. A power substitution extends the classification
to $X^r(X^{d(q-1)}+a)$ whenever $\gcd\bigl(d,(q^e-1)/(q-1)\bigr)=2$.
\end{abstract}
\maketitle

\section{Introduction}\label{sec:intro}

Let $q$ be a prime power. A polynomial over a finite field is a
\emph{permutation polynomial} if its evaluation map is a bijection. The
binomials $X^r(X^{q-1}+a)$ over $\mathbb{F}_{q^e}$ were studied by
Masuda, Rubio and Santiago in several extension degrees, and they
conjectured a general description \cite{MRS22}; a complete
classification of that family is now available \cite{Fan26,NPYY26}. In
this paper we treat the exponent difference $2(q-1)$, that is, the
binomials
\[
 F_{r,a}(X):=X^r\bigl(X^{2(q-1)}+a\bigr)\in\mathbb{F}_{q^e}[X],
 \qquad e,r\geqslant1,
\]
and we determine for all $q$, $e$, $r$ and $a$ whether $F_{r,a}$
permutes $\mathbb{F}_{q^e}$.

The classification for the exponent difference $q-1$ settles this
question whenever $(q^e-1)/(q-1)$ is odd, that is, whenever $q$ is even
or $e$ is odd. The remaining cases have $q$ odd and $e$ even. For the
quadratic extension several families were already known. If
$a^{q+1}=1$, the problem reduces to a monomial on the group of roots of
unity of order dividing $q+1$. The cases $r=1$ and $r=3$ are recorded,
with the earlier sources, in \cite[Results~1.1, 1.2 and 1.4]{HPL23}.
For each fixed odd $r>3$, Hou proved that only finitely many pairs
$(q,a)$ with $q$ odd, $a\in\mathbb{F}_{q^2}^\times$ and $a^{q+1}\ne1$ can
give a permutation polynomial \cite[Theorem~2.1]{Hou16}. Apart from
the case $a^{q+1}=1$, these results concern individual exponents. If $r$ is allowed to vary with $q$, a
classification must determine the residue of $r$ modulo the order
$(q+1)/2$ of the group on which the induced map acts, a modulus that
grows with $q$; this is the main difficulty in the quadratic case.

\subsection{Main results}\label{subsec:results}
For an integer $n\geqslant1$, write
\[
 \ell_n:=\frac{q^n-1}{q-1}.
\]
For a positive integer $h$, the notation $[u]_h$ denotes the least
nonnegative residue of an integer $u$ modulo $h$.

Our first result shows that extensions of even degree at least four
contain no permutation binomials of this form.

\begin{theorem}\label{thm:higher}
Let $q$ be an odd prime power and let $e\geqslant4$ be even. For every
$r\geqslant1$ and every $a\in\mathbb{F}_{q^e}^\times$, the binomial
$F_{r,a}$ does not permute $\mathbb{F}_{q^e}$.
\end{theorem}

In the quadratic extension we obtain an explicit classification in
which the exponent is unrestricted.

\begin{theorem}\label{thm:quadratic}
Let $q$ be an odd prime power, let $r\geqslant1$, and let
$a\in\mathbb{F}_{q^2}^\times$. Put
\[
 h:=\frac{q+1}{2},\qquad s:=(-a)^h,\qquad j:=[r-2]_h.
\]
Then $F_{r,a}$ permutes $\mathbb{F}_{q^2}$ if and only if
$\gcd(r,q-1)=1$ and one of the following conditions holds:
\begin{enumerate}[label=\textup{(\roman*)},leftmargin=*]
\item $s=-1$ and $\gcd(r-2,q+1)=1$;
\item $s^2\ne1$ and $(q,j,s)$ satisfies a row of
Table~\ref{tab:quadratic}.
\end{enumerate}
\end{theorem}

\begin{table}[ht]
\centering
\caption{The quadratic cases with $s^2\ne1$. All conditions in a row
are imposed together.}\label{tab:quadratic}
\begin{tabular}{@{}lll@{}}
\toprule
Restriction on $q$ & Exponent residue $j$ & Coefficient condition\\
\midrule
$q$ odd & $h-1$ & $s=3$\\
$q$ odd & $0$ & $s^2=-1$\\
$q$ odd & $1$ & $3s=1$\\
$q\equiv1\pmod4$ & $(h-1)/2$ & $s^2-2s+5=0$\\
$q\equiv1\pmod4$ & $(h+1)/2$ & $5s^2-2s+1=0$\\
$q=11$ & $3$ & $s\in\{2,6\}$\\
\bottomrule
\end{tabular}
\end{table}

The rows of Table~\ref{tab:quadratic} are alternatives; for small fields
a pair $(j,s)$ may satisfy more than one row. Integers in the
coefficient conditions are interpreted in $\mathbb{F}_{q^2}$, so these
conditions make sense in characteristics three and five without division
by either prime. The condition $\gcd(r,q-1)=1$ is imposed in every row.
In particular, the row with $j=0$ occurs only when $q\equiv1\pmod4$: if
$q\equiv3\pmod4$, then $h$ is even, and $r\equiv2\pmod h$ would make $r$
even, contrary to $\gcd(r,q-1)=1$.

The first and third rows contain the exponents $r=1$ and $r=3$, for
which $r\equiv1$ and $r\equiv3\pmod h$. For a fixed odd $r>3$, the
theorem gives an explicit bound: if $a^{q+1}\ne1$ and $F_{r,a}$ permutes
$\mathbb{F}_{q^2}$, then $q\leqslant4r-7$. Indeed, if $r-2<h$ then
$j=r-2\notin\{0,1\}$, and the remaining rows require $r-2$ to equal
$h-1$, $(h-1)/2$, $(h+1)/2$, or $3$ with $q=11$; in each case
$q\leqslant4r-7$. If $r-2\geqslant h$, then $q\leqslant2r-5$. This is an
explicit form of Hou's finiteness theorem. The last row is sporadic: in
$\mathbb{F}_{121}=\mathbb{F}_{11}[T]/(T^2-2)$ the coefficient $a=5+7T$
has $(-a)^6=2$, and $X^{11}(X^{20}+a)$ permutes $\mathbb{F}_{121}$.

Together with the known results, Theorems~\ref{thm:higher}
and~\ref{thm:quadratic} give the complete classification.

\begin{theorem}\label{thm:classification}
Let $q$ be a prime power, let $e,r$ be positive integers, and let
$a\in\mathbb{F}_{q^e}$.
\begin{enumerate}[label={(\arabic*)},leftmargin=*]
\item\label{class:zero} If $a=0$, then $F_{r,a}$ permutes $\mathbb{F}_{q^e}$ if and
only if $\gcd(r+2(q-1),q^e-1)=1$.
\item\label{class:base} If $a\ne0$ and $e=1$, then $F_{r,a}$ permutes $\mathbb{F}_q$
if and only if $a\ne-1$ and $\gcd(r,q-1)=1$.
\item\label{class:odd-index} Suppose $a\ne0$, $e\geqslant2$, and either $q$ is even or $e$
is odd. Then $F_{r,a}$ permutes $\mathbb{F}_{q^e}$ if and only if
\[
 \gcd(r,q-1)=1,\qquad (-a)^{\ell_e}\ne1,
 \qquad r\ell_u\equiv2\pmod{\ell_e}
\]
for some $u$ satisfying $1\leqslant u<e$ and $\gcd(u,e)=1$.
\item\label{class:higher-even} If $a\ne0$, $q$ is odd, and $e\geqslant4$ is even, then
$F_{r,a}$ does not permute $\mathbb{F}_{q^e}$.
\item\label{class:quadratic} If $a\ne0$, $q$ is odd, and $e=2$, then $F_{r,a}$
permutes $\mathbb{F}_{q^2}$ if and only if the conditions of
Theorem~\ref{thm:quadratic} hold.
\end{enumerate}
\end{theorem}

A monomial substitution (cf.\ \cite[Section~2]{HPL23}) reduces the
exponent differences $d(q-1)$ with $\gcd(d,\ell_e)=2$ to the case $d=2$.
For $r,d\geqslant1$ and $a\in\mathbb{F}_{q^e}$, put
\[
 f_{r,d,a}(X):=X^r\bigl(X^{d(q-1)}+a\bigr).
\]

\begin{corollary}\label{cor:gcd-two}
Let $q$ be a prime power, let $e,r,d\geqslant1$, and let
$a\in\mathbb{F}_{q^e}$, with $\gcd(d,\ell_e)=2$.
Then $q$ is odd and $e$ is even. If $a=0$, the polynomial
$f_{r,d,a}$ permutes $\mathbb{F}_{q^e}$ if and only if
$\gcd(r+d(q-1),q^e-1)=1$.

Suppose $a\ne0$. If $e\geqslant4$, then $f_{r,d,a}$
does not permute $\mathbb{F}_{q^e}$. If $e=2$, let $h=(q+1)/2$ and
$s=(-a)^h$ be as in Theorem~\ref{thm:quadratic}, and put
\[
 \delta:=\frac d2,\qquad j:=[r\delta^{-1}-2]_h,
\]
where $\delta^{-1}$ is taken modulo $h$. Then $f_{r,d,a}$
permutes $\mathbb{F}_{q^2}$ if and only if $\gcd(r,q-1)=1$
and one of the following conditions holds:
\begin{enumerate}[label=\textup{(\roman*)},leftmargin=*]
\item $s=-1$ and $\gcd(r-d,q+1)=1$;
\item $s^2\ne1$ and $(q,j,s)$ satisfies a row of
Table~\ref{tab:quadratic}.
\end{enumerate}
\end{corollary}

\subsection{Ideas of the proofs}\label{subsec:ideas}
For a prime power $M$ and a divisor $D$ of $M-1$, the binomial
$X^r(X^D+a)$ permutes $\mathbb{F}_M$ exactly when $\gcd(r,D)=1$, it has
no nonzero root, and an induced rational function permutes the group of
roots of unity of order $(M-1)/D$ (Lemma~\ref{lem:fibers}). When $q$ is
odd and $e$ is even, a permutation $F_{r,a}$ also forces
$a^{\ell_e/2}\in\mathbb{F}_q^\times$ (Lemma~\ref{lem:coefficient-product});
this restriction on the coefficient is used at every later stage.

\emph{The quadratic case.} For $e=2$ the induced map is
\[
 H(y)=y^j\left(\frac{1+a^qy}{y+a}\right)^2\qquad(y^h=1),
\]
and, once $\gcd(r,q-1)=1$, whether it permutes the roots of unity of
order $h$ depends only on $(j,s)$. We proceed in two stages. First we
determine the coefficient. The power sums $\sum_{y^h=1}H(y)^m$ vanish for
$1\leqslant m<h$, and a residue computation turns them into identities
between $s$ and residues of $j$ (Lemma~\ref{lem:quad-power-sums}). For
$m=1$ the identity reads
\[
 (2j+1)s^2-s-2j+2=0\quad\text{in }\mathbb{F}_q.
\]
After the reciprocity $(j,s)\mapsto([-j]_h,s^{-1})$ we may assume
$0<j\leqslant h/2$ when $j\ne0$. For $0<j<h/2$, combining this identity
with those for $m=2,3,4$ leaves only $3s=1$ and $s^2-2s+5=0$, and the
central residue $j=h/2$ forces $q=11$
(Proposition~\ref{prop:quad-coefficient}). Thus $s$ takes one of a few
values, independently of $q$.

Second we determine $j$; this is the heart of the argument. When
$\alpha$ is a power $K=p^k$ of $p$, or $K+2$, or a small odd multiple of
$K$, the power sum of $F_{r,a}^{(q-1)(q-\alpha)}$ over $\mathbb{F}_{q^2}$
vanishes, and by Lucas' theorem its binomial expansion has very few
nonzero terms. The resulting identity, \eqref{eq:quad-block-one} in the
basic case $\alpha=K$, expresses the quotient on division by $K$ of a
residue determined by $r$ in terms of $s$ and of one bit recording the
size of the remainder. Since $s$ is known,
each identity fixes one more base-$p$ digit of $j$ and bounds what
remains, and an induction on $k$ determines $j$ completely. The
induction depends on the arithmetic of $s$. If $s\notin\mathbb{F}_p$, the
Frobenius map interchanges the two roots of $s^2-2s+5$ and the
identities alternate with the parity of $k$
(Proposition~\ref{prop:quad-nonsplit}). For $s=1/3$ and for split roots
of $s^2-2s+5$, each stage combines several of these identities, and in
the split case also power sums of $H$ of order about $K/2$
(Lemmas~\ref{lem:quad-linear-exponent} and~\ref{lem:quad-split-exponent}).
In characteristic five both coefficient conditions become $s=2$, and both
$j=1$ and $j=(q-1)/4$ survive (Lemma~\ref{lem:quad-five-exponent}).
Sufficiency is verified directly on roots of unity.

\emph{Even degree at least four.} If $4\mid e$, a permutation of
$\mathbb{F}_{q^e}$ yields one of $\mathbb{F}_{q^{e/2}}$
(Lemma~\ref{lem:higher-descent}), and $e=4$ is excluded by comparing two
quadratic coefficients whose parameters $s$ have opposite signs,
together with one further power sum
(Proposition~\ref{prop:higher-quartic}). It remains to treat $e=2u$ with
$u\geqslant3$ odd. Here power sums force, after inversion, a congruence
$r\ell_t\equiv2\pmod{\ell_u}$ with $\gcd(t,u)=1$
(Proposition~\ref{prop:higher-exponent}). The proof turns the vanishing
of suitable power sums into conditions on the base-$q$ digits of a
residue attached to $r$, and a combinatorial argument shows that these
digits form the lower word of a line of rational slope
(Lemma~\ref{lem:higher-digit-word}). Given this congruence, fix a generator $\kappa$ of
$\mathbb{F}_{q^2}^\times$ and $c\in\mathbb{F}_{q^u}^\times$. The equation
$F_{r,a}(x)=F_{r,a}(\kappa cx)$ has a solution
$x\in\mathbb{F}_{q^{2u}}^\times$ exactly when the norm to
$\mathbb{F}_{q^2}$ of $(\lambda-b)/(\lambda-d)$, where $\lambda=c^r$ and
$b,d\in\mathbb{F}_{q^2}$ depend only on $\kappa$ and $r$, lies in a
prescribed class (Lemma~\ref{lem:higher-collision}). We show that this class is attained
by counting the values of the norm ratio on
$\mathbf P^1(\mathbb{F}_{q^u})$ through a weighted Euler product over
$\mathbb{F}_q$ and a transfer matrix on cyclic words
(Proposition~\ref{prop:higher-norm}); for $q=3$ an explicit word suffices
(Corollary~\ref{cor:higher-norm-three}).

\subsection{Organization}
Section~\ref{sec:reductions} contains the reductions and Lucas' theorem.
Section~\ref{sec:quadI} studies the induced map in the quadratic case,
restricts the coefficient and settles the middle residues;
Section~\ref{sec:quadII} determines the exponent and proves
Theorem~\ref{thm:quadratic}. Section~\ref{sec:higher} proves
Theorem~\ref{thm:higher}, and Section~\ref{sec:gcd-two} deduces
Theorem~\ref{thm:classification} and Corollary~\ref{cor:gcd-two}.

\section{Preliminaries}\label{sec:reductions}

For a positive integer $d$ prime to the characteristic, let $\mu_d$
denote the group of $d$-th roots of unity in an algebraic closure.
Whenever $\mu_d$ occurs below, it is contained in the field under
consideration.

\begin{lemma}\label{lem:fibers}
Let $M$ be a prime power, let $D\mid M-1$, and put
$\lambda:=(M-1)/D$. For $a\in\mathbb{F}_M^\times$, the polynomial
$X^r(X^D+a)$ permutes $\mathbb{F}_M$ if and only if
\[
 \gcd(r,D)=1,\qquad (-a)^\lambda\ne1,
\]
and $y\mapsto y^r(y+a)^D$ permutes $\mu_\lambda$.
Coefficients $a,b$ with $a^\lambda=b^\lambda$ give equivalent
permutation conditions for the same $r$.
\end{lemma}
\begin{proof}
The binomial has a nonzero root exactly when $-a$ is a $D$-th power in
$\mathbb{F}_M^\times$, that is, when $(-a)^\lambda=1$. Assume that this
does not happen, so that $\mathbb{F}_M^\times$ is mapped into itself.
On each fiber of $x\mapsto x^D$, replacing $x$ by $vx$ with
$v\in\mu_D$ multiplies the value by $v^r$. Hence the map sends each
fiber bijectively onto a fiber exactly when $\gcd(r,D)=1$. Raising its
value to the $D$-th power gives the stated map on the set of fibers,
which $x\mapsto x^D$ identifies with $\mu_\lambda$.
Finally, $a^\lambda=b^\lambda$ means $a=bc^D$ for some
$c\in\mathbb{F}_M^\times$; substituting $cX$ for $X$ and dividing by
$c^{r+D}$ replaces $a$ by $b$.
\end{proof}

Composition with the inversion of $\mathbb{F}_M^\times$, extended by
$0\mapsto0$, also preserves the permutation property. Since
$x^{-r}(x^{-D}+a)=a\,x^{-r-D}(x^D+a^{-1})$, dividing the output by $a$
replaces the parameters by
\begin{equation}\label{eq:reciprocity}
 (r,a)\longmapsto(r',a^{-1}),\qquad
 r'\equiv-r-D\pmod{M-1},\quad r'\geqslant1.
\end{equation}
Choosing $r'\geqslant1$ makes the new binomial a polynomial that
vanishes at $0$.

\begin{lemma}\label{lem:coefficient-product}
Suppose $q$ is odd, $e$ is even, and $a\in\mathbb{F}_{q^e}^\times$.
Put $D:=2(q-1)$ and $\lambda:=(q^e-1)/D$. If $F_{r,a}$ permutes
$\mathbb{F}_{q^e}$, then
\[
 \gcd(r,D)=1,\qquad a^\lambda\in\mathbb{F}_q^\times.
\]
In particular, $a$ is a square in $\mathbb{F}_{q^e}$.
\end{lemma}
\begin{proof}
The first assertion follows from Lemma~\ref{lem:fibers}; in
particular $r$ is odd. Set $w=a^\lambda$ and
$\epsilon=(-1)^\lambda$. Taking the product of the induced
permutation over $\mu_\lambda$ gives $(w-\epsilon)^D=1$.
Here $\prod_{y\in\mu_\lambda}(y+a)=w-\epsilon$, and the product
of the powers $y^r$ equals the product of the $y$, since $r$ is odd.
We also have $w^D=1$. Consequently the $q$-th power of each of
$w$ and $w-\epsilon$ is itself or its negative.
If $w^q=-w$, neither choice is possible for $w-\epsilon$:
the two equalities would respectively imply $2w=0$ or
$2\epsilon=0$. Thus $w^q=w$. It follows that
$a^{(q^e-1)/2}=w^{q-1}=1$.
\end{proof}

Several arguments use Lucas' theorem in the following form:
\[
 \binom AB=\prod_{i\geqslant0}\binom{A_i}{B_i}
 \quad\text{in }\mathbb{F}_p,
 \qquad A=\sum_i A_ip^i,\quad B=\sum_i B_ip^i,
 \quad 0\leqslant A_i,B_i<p.
\]
Here a factor is zero if $B_i>A_i$. The congruence follows by
comparing coefficients in
$(1+X)^A=\prod_i(1+X^{p^i})^{A_i}$ in characteristic $p$.

\section{The quadratic case: the induced map and the coefficient}\label{sec:quadI}

Throughout this section and the next, $q=p^f$ is odd,
$a\in\mathbb{F}_{q^2}^\times$, and $h$, $s$, $j$ are as in
Theorem~\ref{thm:quadratic}:
\[
 h=\frac{q+1}{2},\qquad s=(-a)^h,\qquad j=[r-2]_h.
\]

\subsection{The induced map and its power sums}\label{sec:quad-quotient}
Since $(q^2-1)/(2(q-1))=h$, Lemma~\ref{lem:coefficient-product} shows
that $s\in\mathbb{F}_q^\times$ whenever $F_{r,a}$ is a permutation; in
all cases $s^2=a^{q+1}$. The absence of
nonzero roots is equivalent to $s\ne1$. Since $q-1$ is even,
$\gcd(r,2(q-1))=1$ is equivalent to $\gcd(r,q-1)=1$.

Using $y^q=y^{-1}$ on $\mu_h$, the map $y\mapsto y^r(y+a)^{2(q-1)}$
of Lemma~\ref{lem:fibers} becomes
\begin{equation}\label{eq:quad-quotient}
 H(y):=y^j\left(\frac{1+a^qy}{y+a}\right)^2,
 \qquad y\in\mu_h.
\end{equation}
Consequently, once $\gcd(r,q-1)=1$ is imposed, the permutation
condition depends only on $(j,s)$. Formula~\eqref{eq:reciprocity}
replaces this pair by
\begin{equation}\label{eq:quad-reciprocity}
 (j,s)\longmapsto([-j]_h,s^{-1}).
\end{equation}

If $H$ permutes $\mu_h$, then $\sum_{y\in\mu_h}H(y)^m=0$ for
$1\leqslant m<h$. The next lemma expresses these power sums through
binomial sums.

\begin{lemma}\label{lem:quad-power-sums}
Suppose $s\in\mathbb{F}_q^\times$, $s^2\ne1$, and $H$ in
\eqref{eq:quad-quotient} permutes $\mu_h$. For $1\leqslant m<h$, put
$d=2m-1$ and write $mj=bh+K$, $0\leqslant K<h$. Set
\[
 \lambda=\frac{s^2-1}{s^2},\qquad
 D_\lambda(d,n)=\sum_{t=0}^d\binom dt\binom nt\lambda^t,
\]
using integer binomial coefficients also when $n=-1$. Then
\begin{align}
 R_{m,K}&:=s^{2d}\bigl(D_\lambda(d,K-1)
                   +sD_\lambda(d,K+h-1)\bigr)
                      \label{eq:quad-binomial-residue}\\
 &=\begin{cases}0,&K>0,\\1,&K=0.\end{cases}
                      \label{eq:quad-residue-value}
\end{align}
Equivalently,
\begin{equation}\label{eq:quad-residue}
 R_{m,K}=\frac{(1-s^2)^{2m}}{2^d(s-1)}
 [z^d](1+z)^{2K-1}
 \left(1-\frac{s}{1-s}z\right)^d
 \left(1+\frac{s}{1+s}z\right)^{2m}
 \left(1+\frac z2\right)^{-2m}.
\end{equation}
\end{lemma}
\begin{proof}
Put $x=-y/a$. Then $x^h=s^{-1}$, and, up to a nonzero constant,
$H(y)=x^j((1-s^2x)/(1-x))^2$. For the differential
\[
 \omega=\frac{x^{K-1}(1-s^2x)^{2m}}
 {(sx^h-1)(1-x)^{2m}}\,dx,
\]
the sum of the residues at $sx^h=1$ is zero by the permutation
hypothesis. The residue at infinity is zero since $K<h$, and
the residue at zero is $-\mathbf1_{K=0}$. Thus
$\operatorname{Res}_{x=1}\omega=\mathbf1_{K=0}$.

At $x=1+T$, multiply by $s(1+T)^h+1$. Since $2h=q+1$,
\[
 s^2(1+T)^{2h}-1\equiv s^2(T+\lambda)\pmod{T^q}.
\]
As $d<q$ and $\lambda\ne0$, cancellation gives
\[
 \operatorname{Res}_{x=1}\omega
 =s^{2d}[T^d](1+T)^{K-1}(T+\lambda)^d
                   \bigl(1+s(1+T)^h\bigr).
\]
Expanding $(T+\lambda)^d$ proves \eqref{eq:quad-binomial-residue}.
Alternatively, substitute $x=(1+z)^2$ and use
$(1+z)^{2h}\equiv1+z\pmod{z^q}$. Factoring
$1-s^2(1+z)^2$ and $1-(1+z)^2=-2z(1+z/2)$ gives
\eqref{eq:quad-residue}.
\end{proof}

\begin{lemma}\label{lem:quad-first}
Suppose $F_{r,a}$ permutes $\mathbb{F}_{q^2}$ and $s^2\ne1$.
If $j=0$, then $s^2=-1$. If $0<j<h$, then
\begin{equation}\label{eq:quad-first}
 (2j+1)s^2-s-2j+2=0.
\end{equation}
\end{lemma}
\begin{proof}
Apply \eqref{eq:quad-binomial-residue} with $m=1$ and
$D_\lambda(1,n)=1+n\lambda$. Since $h=1/2$ in the field,
\[
 R_{1,j}=\frac{s+1}{2}\bigl((2j+1)s^2-s-2j+2\bigr)
 \quad(j>0),\qquad R_{1,0}-1=\frac{s(s^2+1)}2.
\]
Now use \eqref{eq:quad-residue-value} and $s(s+1)\ne0$.
\end{proof}

If $F_{r,a}$ is a permutation and $s^2\ne1$, then
$s\in\mathbb{F}_{p^{\gcd(2,f)}}$. Indeed $s\in\mathbb{F}_q$, and by
Lemma~\ref{lem:quad-first} it is a root of a nonzero polynomial over
$\mathbb{F}_p$ of degree at most two: $X^2+1$ if $j=0$, and the
polynomial in \eqref{eq:quad-first}, whose coefficient of $s$ is $-1$,
if $j>0$.

\subsection{Three elementary families}\label{subsec:quad-elementary}

\begin{lemma}\label{lem:quad-elementary-families}
Assume $\gcd(r,q-1)=1$.
\begin{enumerate}[label={(\arabic*)},leftmargin=*]
\item\label{quad-family:norm-one} If $a^{q+1}=1$, then $F_{r,a}$ permutes $\mathbb{F}_{q^2}$
if and only if $s=-1$ and $\gcd(r-2,q+1)=1$.
\item\label{quad-family:zero} If $s^2\ne1$ and $j=0$, then $F_{r,a}$ permutes
$\mathbb{F}_{q^2}$ if and only if $s^2=-1$.
\item\label{quad-family:linear} The conditions $j=h-1$, $s=3$, and the reciprocal conditions
$j=1$, $3s=1$, are sufficient whenever $s^2\ne1$.
\end{enumerate}
\end{lemma}
\begin{proof}
For~\ref{quad-family:norm-one}, the induced map for $x\mapsto x^{q-1}$ acts on
$\mu_{q+1}$. Once nonzero roots are excluded, the identity
\[
 y^r(y^2+a)^{q-1}=a^{-1}y^{r-2}
\]
proves the assertion. Indeed $s^2=1$, so absence of roots says
$s=-1$; the remaining map is a monomial.

The necessity in~\ref{quad-family:zero} is Lemma~\ref{lem:quad-first}. Conversely,
if $a^{q+1}=-1$ and $j=0$, then
\[
 H(y)=a^{-2}\left(\frac{a-y}{a+y}\right)^2.
\]
The fractional linear map in this expression is injective on
$\mu_h$. Two of its values could be negatives only if their
arguments $y,z$ satisfied $yz=a^2$. This is impossible because
$(yz)^h=1$ whereas $(a^2)^h=-1$. Thus $H$ is injective,
and Lemma~\ref{lem:fibers} applies.

For~\ref{quad-family:linear}, the known criterion for
$X(X^{2(q-1)}+a)$ with $a^{q+1}\ne1$ is $s=3$
\cite[Result~1.2]{HPL23}. Its induced map is
\eqref{eq:quad-quotient} with $j=h-1$, so the same criterion
is sufficient for every $r$ in that residue class satisfying
the fiber coprimality condition. Inversion gives the second case.
\end{proof}

\subsection{Restricting the coefficient}\label{sec:quad-coefficient}

We now use Lemma~\ref{lem:quad-power-sums} for small $m$ to restrict
$s$. By \eqref{eq:quad-reciprocity}, it suffices to treat
$0<j\leqslant h/2$ when $j\ne0$.

\begin{proposition}\label{prop:quad-coefficient}
Assume that $F_{r,a}$ permutes $\mathbb{F}_{q^2}$, that $s^2\ne1$, and,
after applying \eqref{eq:quad-reciprocity} if necessary, that
$0<j\leqslant h/2$. If $0<j<h/2$, then one of the following holds:
\begin{enumerate}[label=\textup{(\roman*)},leftmargin=*]
\item $p\ne3$, $3s=1$, and $j\equiv1\pmod p$;
\item $s^2-2s+5=0$ and $4j\equiv-1\pmod p$.
\end{enumerate}
For $p=3$, the second alternative holds with the stronger congruence
$j\equiv2\pmod9$.
If $j=h/2$, then $q=11$ and $s\in\{2,6\}$.
Conversely, if $q=11$, $j=3$ and $s\in\{2,6\}$, then $F_{r,a}$ permutes
$\mathbb{F}_{121}$ whenever $\gcd(r,10)=1$.
\end{proposition}

\begin{proof}
\emph{The range $0<j<h/2$ when $p\geqslant5$.}
Let $p\geqslant5$ and $0<j<h/2$. Set
\[
 A(s)=s^2-2s+5,\qquad B(s)=s^2+2s-1,\qquad
 J(s)=-\frac{s^2-s+2}{2(s^2-1)}.
\]
By \eqref{eq:quad-first}, $j=J(s)$ in $\mathbb{F}_q$. For $p>2m-1$
and $K>0$, the coefficient in \eqref{eq:quad-residue} may therefore
be calculated with exponent
\[
 2K-1=2mJ(s)-b-1.
\]
Let $P_{m,b}(s)$ be this coefficient multiplied by
$(1-s^2)^{2m-1}$. It is a polynomial: in the term contributing to
$z^{2m-1}$, the denominators in $s$ cancel against this displayed
power. The remaining numerical denominators have prime factors at
most $2m-1$, together with powers of two.
For $m=2$, we have $b=0$, $K>0$, and direct expansion gives
\begin{equation*}
 6P_{2,0}(s)=s(3s-1)A(s)B(s).
\end{equation*}
Under $3s=1$, \eqref{eq:quad-first} gives $j\equiv1\pmod p$.
Under $A(s)=0$, it becomes $(4j+1)(s-3)=0$; since $A(3)=8$,
this gives $4j\equiv-1\pmod p$.

It remains to remove $B(s)=0$ except where it also satisfies $A(s)=0$.
Under $B(s)=0$, \eqref{eq:quad-first} reads $4js=3-3s$. In characteristic
five this forces $s\in\mathbb{F}_5$, whereas $B$ has nonsquare
discriminant $8$ there. Suppose henceforth that $p\geqslant7$.
If $3j<h$, use $(m,b)=(3,0)$; if $3j>h$, use $(m,b)=(3,1)$.
For $m=3,4$, reduction of \eqref{eq:quad-residue} modulo
$B(s)=s^2+2s-1$ gives
\[
\begin{array}{c|c|r}
 (m,b)&P_{m,b}(s)\bmod B(s)&\text{norm of the linear factor}\\ \hline
 (3,0)&\frac35(17-41s)&2\\
 (3,1)&-\frac35(-4317+10426s)&-17\cdot31\cdot89\\
 (4,1)&-\frac{136}{5}(5451-13160s)&-1879.
\end{array}
\]
The entries follow by expanding through degree $2m-1\leqslant7$
and using $s^2=1-2s$. The norm of $u+vs$ modulo $B$ is
$u^2-2uv-v^2$. The first row excludes $B(s)=0$ when $3j<h$.
If $3j>h$, the second row leaves only $p=17,31,89$;
for these primes $p>7$ and $h<4j<2h$, so the third row,
with $(m,b)=(4,1)$, applies. Its norm $-1879$ excludes $31$ and $89$.
At $p=17$ the second row gives
$-4317+10426s\equiv1+5s=11(2s-3)=0$, hence $s=3/2$;
thus $A(s)=0$ by $A-B=6-4s$. If $3j=h$, then \eqref{eq:quad-first} reads
$4s^2-3s+5=0$, whose resultant with $B$ is $158=2\cdot79$.
But $79^f\equiv1\pmod3$, so this characteristic cannot have
$3\mid h$. This settles the range $0<j<h/2$ for $p\geqslant5$.

\emph{The range $0<j<h/2$ when $p=3$.}
In characteristic three, \eqref{eq:quad-first} leaves
\[
 j=3v,\ s^2=s+1,\qquad\text{or}\qquad j=3v+2,\ s^2=1-s.
\]
Both quadratics are irreducible over $\mathbb F_3$, so $f$ is even
and $h\equiv5\pmod9$. For $m=2$, we have $K=2j>0$ and
$D_\lambda(3,n)=1+\lambda^3\lfloor n/3\rfloor$ in
$\mathbb F_9$, by $(1+T)^3=1+T^3$.
The two quadratics give $\lambda^3=-s,s$, respectively.
Substitution in \eqref{eq:quad-binomial-residue} yields
\[
 s^{-6}R_{2,2j}=
 \begin{cases}v+(1-v)s,&j=3v,\ s^2=s+1,\\
              -v,&j=3v+2,\ s^2=1-s.
 \end{cases}
\]
The first value is nonzero because $s\notin\mathbb F_3$.
The second forces $3\mid v$, hence $j\equiv2\pmod9$.

\emph{The central residue $j=h/2$.}
Suppose $j=h/2$. Then $q\equiv3\pmod4$.
The case $q=3$ has no $s\in\mathbb{F}_3^{\times}$ with $s^2\ne1$.
Otherwise \eqref{eq:quad-first} gives
\[
 Q(s)=3s^2-2s+3=0.
\]
In characteristic three, $Q(s)=0$ forces $s=0$. For
$p\geqslant5$, the $m=2$, $K=0$ case of
\eqref{eq:quad-binomial-residue}, using $h-1=-1/2$ in $\mathbb{F}_p$,
gives
\[
 16(R_{2,0}-1)=s(s^2+1)(5s^4-2s^2+5).
\]
As $2s=3(s^2+1)$, we have $s^2+1\ne0$. For $t=s^2\ne0$,
the two equations are $5t^2-2t+5=9t^2+14t+9=0$;
subtracting five times the second from nine times the first
gives $-88t=0$, hence $p=11$.
The roots of $Q$ are then $s=2,6$. Suppose $q>11$. Take $m=6<h$, $K=0$, and $d=11<q$.
As $s\in\mathbb{F}_{11}$, the only nonzero terms in
$D_\lambda(11,n)$ are $1+\binom n{11}\lambda$.
The values $\binom{-1}{11}=-1$ and
$\binom{h-1}{11}\equiv5\pmod{11}$ (by Lucas' theorem) give
\[
 R_{6,0}-1=s^3(1+5\lambda)
             =\frac{s(s^2+1)}2\ne0,
\]
where $s^2+1\ne0$ follows from $Q(s)=0$.
This contradicts \eqref{eq:quad-residue-value}; hence $q=11$.

\emph{The case $q=11$.}
We work in $\mathbb F_{121}=\mathbb F_{11}[T]/(T^2-2)$.
With $w=6+10T$ of order six and $a=5+7T$, one has
$(-a)^6=2$; for $j=3$ the exponents of $H(w^k)$ relative to
$w$ are $(4,5,3,2,0,1)$ for $k=0,\ldots,5$.
Thus $H$ permutes $\mu_6$. The reciprocity
\eqref{eq:quad-reciprocity}, which sends $(j,s)$ to
$([-j]_6,s^{-1})$, also covers $s=6=2^{-1}$, since
$[-3]_6=3$. By the last assertion of Lemma~\ref{lem:fibers}, the
permutation property depends on $a$ only through $(-a)^6=s$, and the
fiber condition $\gcd(r,20)=1$ is equivalent to $\gcd(r,10)=1$. This
completes the proof.
\end{proof}

\subsection{The middle residues}\label{subsec:quad-middle}

\begin{proposition}\label{prop:quad-middle}
Suppose $q\equiv1\pmod4$, $\gcd(r,q-1)=1$, and $s^2\ne1$.
If $j=(h-1)/2$, then $F_{r,a}$ permutes $\mathbb{F}_{q^2}$ if and
only if $s^2-2s+5=0$. If $j=(h+1)/2$, the corresponding condition
is $5s^2-2s+1=0$.
\end{proposition}

\begin{proof}
For the first index, necessity follows immediately from
\eqref{eq:quad-first}. Suppose $s^2-2s+5=0$.
Its roots belong to $\mathbb{F}_q$, since $-1$ is a square there;
the root zero in characteristic five is excluded by $s\ne0$.
Put $U=\mu_{q+1}$ and
\[
 g(x)=\frac{x(x^4+a)}{a^qx^4+1},\qquad x\in U.
\]
Its numerator and denominator are nonzero on $U$, and
$g(x)^q=g(x)^{-1}$, so $g$ maps $U$ into itself. Moreover,
$g(-x)=-g(x)$. Since $h$ is odd, fourth powering identifies
$U/\{1,-1\}$ with $\mu_h$. Using $4j=q-1$ gives
\[
 H(x^4)=g(x)^{-2}.
\]
It is therefore enough to prove that $g$ is injective.

Suppose $g(tx)=g(x)$ for $x,t\in U$ and $t\ne1$. Canceling
$t-1$ in the cross-multiplied equality gives, with $z=x^4$,
\[
 a^qt^4z^2+
 \bigl(t^4+t^3+t^2+t+1-s^2t(t^2+t+1)\bigr)z+a=0.
\]
Set
\[
 T=t+t^{-1}\in\mathbb{F}_q,\qquad v=\frac{a^qt^2x^4}{s}\in U,
 \qquad D(T)=T^2+(1-s^2)T-(1+s^2).
\]
Dividing the preceding equality by $t^2z$ gives
\[
 s(v+v^{-1})=-D(T).
\]
Write $\chi$ for the quadratic character of $\mathbb{F}_q$.
Since $h$ is odd,
\[
 v^h=-\chi(s),\qquad
 D(T)-2s=\left(T+\frac{1-s^2}{2}\right)^2.
\]
For $v\in U\setminus\{-1\}$,
\[
 \chi(v+v^{-1}+2)
 =\left(\frac{(v+1)^2}{v}\right)^{(q-1)/2}=v^h.
\]
On the other hand,
\[
 s(v+v^{-1}+2)=-\left(T+\frac{1-s^2}{2}\right)^2.
\]
Both sides are nonzero, and $\chi(-1)=1$, so this gives
$\chi(v+v^{-1}+2)=\chi(s)$, a contradiction.

If $v=-1$, then $\chi(s)=1$ and $T=(s^2-1)/2=s-3$.
Consequently $T^2-4=-4s$ is a nonzero square in $\mathbb{F}_q$.
The roots $t,t^{-1}$ of $X^2-TX+1$ are therefore distinct and lie
in $\mathbb{F}_q$. This is impossible because
$U\cap\mathbb{F}_q=\{1,-1\}$, and either of these elements equals
its inverse. Thus $g$ is injective.
The second index follows by \eqref{eq:quad-reciprocity}.
\end{proof}

\section{The quadratic case: determining the exponent}\label{sec:quadII}

We now determine $j$. Suppose that $F_{r,a}$ permutes $\mathbb{F}_{q^2}$,
that $s^2\ne1$ and that $j\ne0$. By Proposition~\ref{prop:quad-coefficient},
after the reciprocity \eqref{eq:quad-reciprocity} we have
$0<j\leqslant h/2$, and either $j=h/2$ and $q=11$, or $0<j<h/2$ and $s$
satisfies $3s=1$ (with $p\ne3$) or $s^2-2s+5=0$.
Section~\ref{sec:quad-sparse} derives the identities used to determine
$j$. Section~\ref{sec:quad-nonsplit} treats $s^2-2s+5=0$ with
$p\equiv3\pmod4$, and Sections~\ref{sec:quad-exponents}--\ref{subsec:quad-five}
treat $s=1/3$ with $p\geqslant7$, $s^2-2s+5=0$ with $p\equiv1\pmod4$ and
$p\geqslant13$, and $s=2$ with $p=5$. Section~\ref{subsec:quad-proof}
completes the proof of Theorem~\ref{thm:quadratic}.

\subsection{Power sums with sparse binomial expansions}\label{sec:quad-sparse}

We keep the notation of Section~\ref{sec:quadI} and assume $s^2\ne1$.
Here we derive identities from power sums of $F_{r,a}$ whose binomial
expansions, by Lucas' theorem, have few nonzero terms. All integer
coefficients in these identities are reduced modulo $p$. In the
congruences for summation indices only the residue of $r$ modulo $h$
enters, so $r$ may be replaced by $j+2$ there.

\begin{lemma}\label{lem:quad-sparse}
Suppose $F_{r,a}$ permutes $\mathbb{F}_{q^2}$. For an odd integer
$\alpha$ with $0<\alpha<q-1$, put $\beta:=q-1-\alpha$. Then
\begin{equation}\label{eq:quad-hermite-row}
 \sum_{\substack{0\leqslant i\leqslant\alpha,\ 0\leqslant l\leqslant\beta\\
        2(i-l)\equiv(\alpha+1)r\ (\mathrm{mod}\ q+1)}}
 \binom{\alpha}{i}\binom{\beta}{l}a^{-i-ql}=0.
\end{equation}
Let $K=p^k$ with $1\leqslant k<f$, and set
\[
 b:=\left[-\frac{(K+1)r}{2}\right]_h,\qquad
 b':=[b+K]_h,\qquad \delta:=\frac{b+K-b'}h.
\]
For $0\leqslant n<h$, put
\[
 D_K(n):=1+\left\lfloor\frac nK\right\rfloor+
 s^{-1}\left(1+\left\lfloor\frac{n+h}{K}\right\rfloor
                 -\mathbf{1}_{n=h-1}\right).
\]
Then
\begin{equation}\label{eq:quad-full-block}
 s^{2K}D_K(b)=s^\delta D_K(b').
\end{equation}
If $0\leqslant b\leqslant h-K-2$, write $b=dK+x$ with
$0\leqslant x<K$, and let $\varepsilon$ be $0$ or $1$ according as
$x<(K-1)/2$ or $x\geqslant(K-1)/2$. In that case
\begin{equation}\label{eq:quad-block-one}
 d+1=\frac1{s^{2K}-1}+\frac{1-2\varepsilon}{2(s+1)}.
\end{equation}
\end{lemma}
\begin{proof}
The exponent $\alpha+q\beta=(q-1)(q-\alpha)$ lies strictly
between $0$ and $q^2-1$. The corresponding power sum of a permutation
vanishes. Expanding the two Frobenius factors, a term survives
summation over $\mathbb{F}_{q^2}^\times$ exactly when its exponent
is divisible by $q^2-1$. Dividing that exponent by $q-1$ gives
the congruence in \eqref{eq:quad-hermite-row}; the common nonzero
factor $a^{\alpha+q\beta}$ may be removed.

For $\alpha=K$, only $i=0,K$ have nonzero binomial coefficients.
For $0\leqslant l<q$, the identity
\[
 (1+X)^{q-1-K}\equiv
 \frac1{(1+X)(1+X^K)}\pmod{X^q}
\]
gives
\[
 \binom{q-1-K}{l}=(-1)^l
       \left(1+\left\lfloor\frac lK\right\rfloor\right).
\]
The coefficient is understood to be zero outside its original range.
At $l=q$ it is also zero; this accounts for the subtraction when
$n=h-1$ in $D_K(n)$.
Put $B=-a^{-q}$. Then $B^h=s^{-1}$ and
$a^{-K}B^K=-s^{-2K}$. For $i=0$ the indices $l$ that occur are $b$ and $b+h$;
for $i=K$ they are $b'$ and $b'+h$. Substitution in
\eqref{eq:quad-hermite-row}, followed by division by $B^b$,
gives \eqref{eq:quad-full-block}.

In the stated range there is no reduction of $b+K$ modulo $h$.
Writing
\[
 h=\frac{q/K-1}{2}K+\frac{K+1}{2}
\]
shows that
\[
 D_K(b)=(1+s^{-1})(d+1)+s^{-1}(-1/2+\varepsilon),\qquad
 D_K(b+K)-D_K(b)=1+s^{-1}.
\]
Equation~\eqref{eq:quad-full-block} now gives
\eqref{eq:quad-block-one}. The denominators are nonzero since
$s^2\ne1$.
\end{proof}

In the applications of \eqref{eq:quad-block-one} below, $b$ is first
given by an integer representative modulo $h$; we then check that
$0\leqslant b\leqslant h-K-2$ and write $b=dK+x$ with $0\leqslant x<K$
and $\varepsilon$ as in the lemma. The next lemma treats two residues
outside this range.

\begin{lemma}\label{lem:quad-boundary-residues}
Suppose $s^K=s$. If the residue $b$ in
Lemma~\ref{lem:quad-sparse} is $h-1$, then
$s^2-2s-3=0$. If $b=h-K-1$, then $s^2=-1$.
\end{lemma}
\begin{proof}
For $b=h-1$ the full expression gives
\[
 b'=K-1,\quad \delta=1,\quad
 D_K(b)=\frac12,\quad D_K(b')=1+\frac{3}{2s}.
\]
For $b=h-K-1$ it gives
\[
 b'=h-1,\quad \delta=0,\quad
 D_K(b)=-\frac12,\quad D_K(b')=\frac12.
\]
Substitution in \eqref{eq:quad-full-block}, using $s^{2K}=s^2$,
proves the two assertions.
\end{proof}

\begin{lemma}\label{lem:quad-six-terms}
Assume $s^K=s$, where $K=p^k$ and $1\leqslant k<f$. Put
\[
 c:=\left[-\frac{(K+3)r}{2}\right]_h.
\]
Suppose $c+K+2\leqslant h-2$, and write
$c=dK+x$, $c+h=d_2K+x_2$, with
$0\leqslant x,x_2\leqslant K-3$. Set
\[
 z:=-s^{-2},\quad
 B_z(x):=\sum_{v=0}^2\binom2v z^v\binom{x+v+2}{2},\quad
 W_z(d):=d+1+z(d+2).
\]
If $F_{r,a}$ is a permutation, then
\begin{equation}\label{eq:quad-block-two}
 W_z(d)B_z(x)+s^{-1}W_z(d_2)B_z(x_2)=0.
\end{equation}
\end{lemma}
\begin{proof}
Use $\alpha=K+2$ in \eqref{eq:quad-hermite-row}. By Lucas' theorem,
the indices $i$ with $\binom{K+2}{i}\ne0$ in $\mathbb{F}_p$ have the form
$i=uK+v$, where $u\in\{0,1\}$ and $0\leqslant v\leqslant2$, and then
$\binom{K+2}{i}=\binom2v$. For the coefficients $\binom{q-K-3}{l}$,
expand
\[
 (1+X)^{q-K-3}\equiv
 (1+X)^{-3}(1+X^K)^{-1}\pmod{X^q}.
\]
After the sign $(-1)^l$ is removed, the coefficient at
$l=dK+x$ equals $(d+1)\binom{x+2}{2}$.
The inequalities on $x,x_2$ ensure that adding $v$ causes no
carry, and the bound on $c$ ensures that no index crosses $h$.
The factors contributed by the powers of $a$ are $z^{u+v}$
and $s^{-1}z^{u+v}$, since $s^K=s$. Summing the six terms
for each of the two values of $l$ gives
\eqref{eq:quad-block-two}.
\end{proof}

\begin{lemma}\label{lem:quad-divisibility}
Let $q=p^f=KP$, where $K=p^k\geqslant7$, $P\geqslant K$,
and suppose that $0<\Delta=mK\leqslant K(P-1)/(K-1)$.
Suppose $F_{r,a}$ permutes $\mathbb F_{q^2}$ and one of the
following holds:
\begin{enumerate}[label=\textup{(\roman*)},leftmargin=*]
\item\label{lift:linear} $p\geqslant7$, $s=1/3$, $j=1+\Delta<h/2$;
\item\label{lift:quadratic} $K\equiv P\equiv1\pmod4$, $p\ne5$,
 $s^2-2s+5=0$, $s^K=s$,
 $j=(q-1)/4-\Delta\in(0,h/2)$.
\end{enumerate}
Then $p\mid m$, except possibly in case~\ref{lift:quadratic} when
$p=3$, $m\equiv1\pmod3$, and the integer
\[
 B:=\frac{P-5}{4}K+\frac{K-1}{2}+\frac{mK(K+1)}2
\]
satisfies $B\geqslant h$.
\end{lemma}
\begin{proof}
Since $K\geqslant7$, the size hypothesis gives
$\Delta=mK\leqslant K(P-1)/(K-1)\leqslant7(P-1)/6$.
In case~\ref{lift:linear}, take
\[
 b_*=h-\frac{(3+\Delta)(K+1)}2,
 \qquad c_*=h-\frac{(3+\Delta)(K+3)}2.
\]
Then $-h\leqslant c_*<b_*\leqslant h-K-2$: the only
non-immediate inequality follows from
\[
 c_*+h\geqslant\frac{(5K-21)P-11K-21}{12}
 \geqslant\frac{(K-7)(5K+3)}{12}\geqslant0.
\]
If $b_*\geqslant0$, then $b=b_*$ has remainder $x=K-1$ and quotient
$d=(P-5)/2-m(K+1)/2$ on division by $K$, and
\eqref{eq:quad-block-one}, with $s=1/3$, gives $p\mid m$.
Otherwise $b=h+b_*$. Since $h-b\equiv1\pmod K$, the bound
$b\leqslant h-K-2$ can fail only for $b=h-1$ or $b=h-K-1$, and
Lemma~\ref{lem:quad-boundary-residues} excludes both values.
Thus \eqref{eq:quad-block-one} applies with
$d=P-2-m(K+1)/2$ and $x=(K-1)/2$, and gives
$m\equiv1\pmod p$. Now $c=h+c_*$ lies in $[0,h-K-4]$, and in the
notation of Lemma~\ref{lem:quad-six-terms} we have $x=(K-7)/2$,
$x_2=K-3$, $d\equiv-7/2$ and $d_2\equiv-4\pmod p$.
For $z=-9$, \eqref{eq:quad-block-two} gives
\[
 W_z(-7/2)B_z(-7/2)+3W_z(-4)B_z(-3)=-120\ne0.
\]
Thus $p\mid m$ in case~\ref{lift:linear}.

In case~\ref{lift:quadratic}, the integer $B$ of the statement
represents the residue $b$ of Lemma~\ref{lem:quad-sparse} for this $K$.
It is positive, and the uniform size bound gives
\[
 2h-K-2-B\geqslant\frac{(2K-7)P+4K+1}{12}>0.
\]
If $h-K-1\leqslant B<h$, then $h-B\equiv1\pmod K$
forces $h-B=1$ or $K+1$. Lemma~\ref{lem:quad-boundary-residues}
would give $(s+1)(s-3)=0$ or $s^2=-1$, neither compatible
with $s^2-2s+5=0$ in characteristic $p\ne2,5$.
Thus $b=[B]_h$ satisfies $0\leqslant b\leqslant h-K-2$.
Its remainder modulo $K$ is $(K-1)/2$ if $B<h$ and $K-1$ otherwise,
so $\varepsilon=1$. Since
$1/(s^2-1)-1/[2(s+1)]=-1/4$, equation~\eqref{eq:quad-block-one} requires
$d\equiv-5/4\pmod p$. The quotient $d$ is
$(P-5)/4+m(K+1)/2$, diminished by $(P+1)/2$ if $B\geqslant h$.
Consequently
\[
 B<h\ \Longrightarrow\ m\equiv0\pmod p,
 \qquad B\geqslant h\ \Longrightarrow\ m\equiv1\pmod p.
\]
In the latter case we use Lemma~\ref{lem:quad-six-terms}; its residue
$c$ has the positive representative
\[
 C=\frac{(K+3)(2mK-3)}4<2h-K-4.
\]
Indeed, the same bound yields
\[
 2h-K-4-C\geqslant\frac{(5K-21)P+4K+12}{12}>0.
\]
The interval $h-K-3\leqslant C<h$ cannot occur:
$h-C\equiv(K+11)/4\pmod K$ would force
$h-C=(K+11)/4$, hence $m(K+3)=P+1$,
contrary to $m\equiv1\pmod p$.
Thus $c=[C]_h$ satisfies the hypotheses of
Lemma~\ref{lem:quad-six-terms}, with the following values of $x$, $x_2$
and of $d$, $d_2$ modulo $p$:
\[
\begin{array}{c|cc|cc}
 &x&x_2&d\pmod p&d_2\pmod p\\ \hline
 C<h&(K-9)/4&(3K-7)/4&1/2&0\\
 C\geqslant h&(3K-11)/4&(K-9)/4&0&1/2.
\end{array}
\]
In the two cases the left side of \eqref{eq:quad-block-two} equals,
respectively,
\[
\begin{aligned}
 E_0(s)&=W_z(1/2)B_z(-9/4)+s^{-1}W_z(0)B_z(-7/4)
       =\frac{12(249+73s)}{5^6},\\
 E_1(s)&=W_z(0)B_z(-11/4)+s^{-1}W_z(1/2)B_z(-9/4)
       =\frac{24(527+29s)}{5^6}.
\end{aligned}
\]
Here $s^{-1}=(2-s)/5$, $z=(1+2s)/25$, and $s^2=2s-5$;
these identities give the displayed simplifications directly.
The products of $E_0,E_1$ with their conjugates under
$s\mapsto2-s$ are $1152/5^6$ and $11520/5^6$.
Both are nonzero for $p\geqslant7$, excluding $B\geqslant h$
in that case. For $p=3$ only the stated residual case remains.
\end{proof}

The argument of Lemma~\ref{lem:quad-sparse} also applies to the
exponents $\alpha=tK$ with $t$ odd and small.

\begin{lemma}\label{lem:quad-large-block}
Suppose $F_{r,a}$ permutes $\mathbb{F}_{q^2}$. Let $K=p^k$ with
$1\leqslant k<f$, suppose $s^K=s$, and let $t<p$ be odd. Assume that
\[
 b:=\left[-\frac{(tK+1)r}{2}\right]_h
\]
satisfies $0\leqslant b\leqslant h-tK-2$, and write $b=dK+x$ with
$0\leqslant x<K$. Let $\varepsilon$ be as in Lemma~\ref{lem:quad-sparse},
put $z=-s^{-2}$, and set
\[
 H_t(d):=\sum_{u=0}^t\binom tu z^u\binom{d+u+t}{t}.
\]
Then
\begin{equation}\label{eq:quad-large-block}
 H_t(d)+s^{-1}H_t(d-1/2+\varepsilon)=0.
\end{equation}
\end{lemma}
\begin{proof}
Use $\alpha=tK$ in \eqref{eq:quad-hermite-row}. By Lucas' theorem the
only indices $i$ with $\binom{tK}{i}\ne0$ in $\mathbb{F}_p$ are $i=uK$,
$0\leqslant u\leqslant t$. After the sign $(-1)^l$ is removed, the
coefficient $\binom{q-1-tK}{l}$ equals
$\binom{\lfloor l/K\rfloor+t}{t}$, as follows from
$(1+X)^{q-1-tK}\equiv(1+X)^{-1}(1+X^K)^{-t}\pmod{X^q}$.
For each $i$, the second value of $l$ exceeds the first by $h$, and its
quotient by $K$ exceeds that of the first by $(q/K-1)/2+\varepsilon$,
whose residue modulo $p$ is $-1/2+\varepsilon$. This proves
\eqref{eq:quad-large-block}.
\end{proof}

\subsection{Coefficients outside the prime field}\label{sec:quad-nonsplit}

In this subsection $p\equiv3\pmod4$ and $s^2-2s+5=0$. The discriminant of
this polynomial is $-16$, so $s\notin\mathbb{F}_p$, and the Frobenius
map $u\mapsto u^p$ interchanges $s$ and $2-s$. The identities of
Section~\ref{sec:quad-sparse} therefore take different forms according
to the parity of $k$ in $K=p^k$.

\begin{proposition}\label{prop:quad-nonsplit}
Suppose $p\equiv3\pmod4$, $s^2-2s+5=0$, and $0<j<h/2$.
If $F_{r,a}$ permutes $\mathbb{F}_{q^2}$, then $f$ is even and
$j=(q-1)/4$. Conversely, if $f$ is even, $j=(q-1)/4$ and
$s^2-2s+5=0$, then $F_{r,a}$ permutes $\mathbb{F}_{q^2}$ whenever
$\gcd(r,q-1)=1$.
\end{proposition}
\begin{proof}
Since $s^2-2s+5$ has discriminant $-16$, its roots lie in
$\mathbb{F}_{p^2}\setminus\mathbb{F}_p$, so $f$ is even. Put
$j_0=(q-1)/4$ and $\Delta=j_0-j\geqslant0$. Equation~\eqref{eq:quad-first} implies $p\mid\Delta$. If $p=3$, Proposition~\ref{prop:quad-coefficient}
strengthens this to $9\mid\Delta$; the case $f=2$ then follows
immediately from $0<j<h/2$.

We now treat $p\geqslant7$ and the case $p=3$, $f\geqslant4$,
together. We apply \eqref{eq:quad-block-one} both to $(j,s)$ and, via
\eqref{eq:quad-reciprocity}, to $([-j]_h,s^{-1})$. If $u$ is either
of $s$, $s^{-1}$, its $p$-th power is
$t=(3u-1)/(u-3)$, and
\[
 \frac1{u^2-1}-\frac1{t^2-1}
 =\frac1{2(u+1)}-\frac1{2(t+1)}\ne0.
\]
Indeed, if $B(u)=1/[2(u+1)]$, the difference between
\eqref{eq:quad-block-one} and its image under the Frobenius map is $-2\varepsilon(B(u)-B(t))$ for odd $k$
and $2(1-\varepsilon)(B(u)-B(t))$ for even $k$.
Thus, whenever \eqref{eq:quad-block-one} applies with $K=p^k$, the
quantities $\varepsilon$ and $d$ satisfy
\begin{equation}\label{eq:quad-nonsplit-digits}
\begin{array}{c|ccc}
 &\varepsilon&d+1\ (u=s)&d+1\ (u=s^{-1})\\ \hline
 k\text{ odd}&0&0&-1/2\\
 k\text{ even}&1&-1/4&-5/4.
\end{array}
\end{equation}
Because $f$ is even, the exponents of $p$ in $K=p^k$ and in $P=q/K$
have the same parity; hence this table also applies with $P$ in place
of $K$.

Starting at $k=1$ for $p\geqslant7$ and at $k=2$ for $p=3$,
we prove for $k\leqslant f/2$ that
\begin{equation}\label{eq:quad-nonsplit-induction}
 p^{k+1}\mid\Delta,\qquad
 \Delta\leqslant
 \begin{cases}
 K(P-3)/(2(K-1)),&k\text{ odd},\\
 K(P-1)/(K-1),&k\text{ even},
 \end{cases}
 \qquad K=p^k,\quad P=q/K.
\end{equation}
If $\Delta=0$ there is nothing to prove. Otherwise write
$\Delta=mK$, $m\geqslant1$, and $T=m(K-1)/2$. The initial
bounds $\Delta<q/4$ and the preceding stage of
\eqref{eq:quad-nonsplit-induction} yield the following strict upper
bounds on $T/P$:
\[
\begin{array}{c|cc|c}
 &\text{first stage}&\text{later stages}&\text{parity of }k\\ \hline
 p\geqslant7&p/8&49p/96&\text{odd}\\
 p\geqslant7&\text{---}&7p/24&\text{even}\\
 p=3&\text{---}&2&\text{odd}\\
 p=3&9/8&1&\text{even}.
\end{array}
\]
For later stages use $K/p\geqslant p^2$ when $k$ is odd
and $K/p\geqslant p$ when $k$ is even; the $p=3$ induction
starts at $k=2$. These bounds include $K=P=9$.

\emph{The size bound.}
Apply \eqref{eq:quad-block-one} with $P$ in place of $K$. An integer
representative of the corresponding residue is
\begin{equation*}
 b_*=C-\frac{3(P+1)}4+T,
 \qquad C=\begin{cases}0,&k\text{ odd},\\ h/2,&k\text{ even}.
 \end{cases}
\end{equation*}
When $k$ is odd, if $b_*\geqslant0$ then
$0\leqslant b_*\leqslant h-P-2$ and
$0\leqslant\lfloor b_*/P\rfloor<p-1$ (for $p=3$ this
quotient is at most one). This contradicts the condition
$d\equiv-1\pmod p$ for odd $k$ and $u=s$ in \eqref{eq:quad-nonsplit-digits}.
Thus $b_*=-L$ with $1\leqslant L\leqslant3(P+1)/4$.
For the reciprocal pair $([-j]_h,s^{-1})$ the corresponding
representative is
\[
 b'=h-2(P+1)+L
    =\frac{K-3}{2}P+L-\frac{P+3}{2};
\]
it satisfies $0\leqslant b'\leqslant h-P-2$, and its quotient by $P$ is
either $(K-5)/2$ or $(K-3)/2$. The condition $d'\equiv-3/2\pmod p$ for
odd $k$ and $u=s^{-1}$ selects the latter. Hence $L\geqslant(P+3)/2$, so
$T\leqslant(P-3)/4$, as required.

When $k$ is even, the representative is
\[
 b_*=\frac{K-5}{4}P+\frac{P-1}{2}+T.
\]
It satisfies $0\leqslant b_*\leqslant h-P-2$, since
$h-P-2-b_*=(K-1)P/4-1-T>0$ by
the displayed bounds. Moreover,
$\lfloor b_*/P\rfloor=(K-5)/4+v$ with $0\leqslant v<p$;
for $p=3$, $v\leqslant1$ even at $K=9$. The condition
$d\equiv-5/4\pmod p$ for even $k$ and $u=s$
forces $v=0$, and thus $T\leqslant(P-1)/2$.
This completes the corresponding size bound in
\eqref{eq:quad-nonsplit-induction}.

\emph{The divisibility.}
For odd $k$, apply \eqref{eq:quad-block-one} with $K$ itself; the
corresponding residue has representative
\begin{equation*}
 b_*=-\frac{3(K+1)}4+\frac{mK(K+1)}2.
\end{equation*}
By the odd size bound,
$0\leqslant b_*<q(K+1)/(4(K-1))\leqslant q/3
\leqslant h-K-1$, and its quotient by $K$ is
$d=m(K+1)/2-1$. The condition $d+1\equiv0\pmod p$ for odd $k$ and $u=s$
forces $p\mid m$. For even $k$,
Lemma~\ref{lem:quad-divisibility}~\ref{lift:quadratic} yields $p\mid m$
unless $p=3$, $m\equiv1\pmod3$, and the integer $B$ of that lemma
satisfies $B\geqslant h$. We exclude this last case next.

Assume therefore that $p=3$, $m\equiv1\pmod3$ and $B\geqslant h$. Then
$k$ is even, $K\geqslant9$, $s^2+s-1=0$, and $K\equiv P\equiv1\pmod8$.
Take $\alpha=K+6$, $z=-s^{-2}=s^2$, and put
\[
 c=\frac{(K+7)(2mK-3)}4.
\]
The even size bound gives $0<c<2h-K-8$: the numerator
of the upper-bound difference after multiplication by $4(K-1)$
is $2K(K-9)P+(K+7)(K+1)>0$.
The interval $h-K-7\leqslant c<h$ is also impossible.
Indeed, $h-c\equiv(K+23)/4\pmod K$ would force
$m(K+7)=P+1$, contrary to reduction modulo eight.

Lucas' theorem restricts the indices $i$ to
$i=uK+v$, where $u\in\{0,1\}$ and $v\in\{0,3,6\}$.
For $l=dK+x<q$ and $0\leqslant x<K$, the coefficient
$\binom{q-1-K-6}{l}$ satisfies
\[
 (-1)^l\binom{q-1-K-6}{l}
  =(d+1)\binom{x+6}{6},\qquad
 \binom{x+6}{6}=\mathbf1_{x\bmod9\in\{0,1,2\}}
 \quad\text{in }\mathbb F_3,
\]
by $(1+X)^{q-1-K-6}\equiv
(1+X)^{-7}(1+X^K)^{-1}\pmod{X^q}$.
Since $c\equiv6$ and $h\equiv5\pmod9$, only $v=3$
contributes without a shift by $h$, and only $v=0$ with a shift.
Writing $M=m(K+7)/2\equiv2\pmod3$, the relevant decompositions are
\[
\begin{aligned}
 c+3&=(M-1)K+(K-9)/4,\\
 c+h&=(M+(P-1)/2-1)K+(3K-19)/4,\\
 c-h&=(M-(P-1)/2-2)K+(3K-23)/4.
\end{aligned}
\]
The remainders lie in $[0,K)$ even for $K=9$, and
the quotients are $1,2,2$ modulo three.
Consequently the factor $d+1$ leaves just two contributions to
\eqref{eq:quad-hermite-row}, whose sum, up to a nonzero factor, is
\[
 z^3+z/s=s^3\ne0\quad(c<h),\qquad
 z^3+zs=1\ne0\quad(c\geqslant h).
\]
This excludes the last case and completes the induction.
At $k=f/2$ the bounds give $\Delta\leqslant K$,
whereas $pK\mid\Delta$, so $\Delta=0$.
The asserted sufficiency is Proposition~\ref{prop:quad-middle}.
\end{proof}

\subsection{The coefficient \texorpdfstring{$s=1/3$}{s=1/3}}\label{sec:quad-exponents}

\begin{lemma}\label{lem:quad-linear-exponent}
Suppose $q=p^f$, $p\geqslant7$, $s=1/3$, and
$0<j<h/2$. If $F_{r,a}$ permutes $\mathbb{F}_{q^2}$, then $j=1$.
\end{lemma}
\begin{proof}
Put $\Delta=j-1$. Equation~\eqref{eq:quad-first} gives
$p\mid\Delta$, and
\[
 0\leqslant\Delta<q/4.
\]
For $f=1$ this proves the assertion. We may assume $f\geqslant2$ and $\Delta>0$.

\emph{Induction.}
We shall prove successively, for $1\leqslant k\leqslant
\lfloor f/2\rfloor$, that
\begin{equation}\label{eq:quad-linear-induction}
 p^{k+1}\mid\Delta,\qquad
 \Delta\leqslant\frac{p^k(p^{f-k}-1)}{p^k-1}.
\end{equation}
At the beginning of a step write
\[
 K=p^k,\qquad P=q/K\geqslant K,\qquad \Delta=mK,\qquad
 T=\frac{m(K-1)}2.
\]
Here $m$ is a positive integer; the required divisibility by $K$
comes from \eqref{eq:quad-first} when $k=1$ and from the preceding
step otherwise.

First use \eqref{eq:quad-block-one} with $P$ in place of $K$.
An integer representative of the corresponding residue is
\begin{equation*}
 b=h-\frac{3(P+1)}2-T
   =\frac{K-5}{2}P+P-1-T.
\end{equation*}
Indeed $q\equiv-1\pmod h$ and $r\equiv3+\Delta\pmod h$.
At the first step $T<pP/8$. Subsequently the preceding bound in
\eqref{eq:quad-linear-induction} and $K\geqslant p^2$ give
\[
 T<\frac{p(K-1)P}{2(K-p)}\leqslant\frac{p+1}{2}P.
\]
These inequalities show $b\geqslant0$; the displayed expression
also gives $b+P\leqslant h-2$. Write $T=vP+\rho$, where
$0\leqslant\rho<P$. Then
\[
 b=dP+x,\qquad d=\frac{K-5}{2}-v,\qquad x=P-1-\rho,\qquad 0\leqslant v<p.
\]
For $s=1/3$ equation~\eqref{eq:quad-block-one} becomes
\begin{equation}\label{eq:quad-linear-digit}
 d+1=
 \begin{cases}
  -3/4,&x<(P-1)/2,\\
  -3/2,&x\geqslant(P-1)/2.
 \end{cases}
\end{equation}
All equalities between these rational constants and the quotient $d$
are in $\mathbb{F}_p$.

Suppose $x<(P-1)/2$. Then $v\equiv-3/4\pmod p$.
The bound $v<(p+1)/2$ excludes $p\equiv1\pmod4$; when
$p\equiv3\pmod4$, its least representative $(p-3)/4$
exceeds $p/8$, excluding $k=1$. For $k\geqslant2$ we use
Lemma~\ref{lem:quad-large-block} with $t=3$ and $P$ in place of $K$; the
corresponding residue has representative $b_3=b-3P+m$.
The preceding bound on $T$, together with $K\geqslant p^2$,
gives
\[
 0<b_3,\qquad b_3+3P=b+m\leqslant h-2,
 \qquad \frac mP<\frac p{K-p}\leqslant\frac1{p-1}.
\]
Since $x<(P-1)/2$, the bound $m/P<1/(p-1)$ also gives
$x+m<P$. Thus $b_3$ has quotient $d-3$ by $P$, where
$d-3\equiv-19/4\pmod p$, and
\eqref{eq:quad-large-block} gives, for some
$\varepsilon_3\in\{0,1\}$,
\[
 H_3(-19/4)+3H_3(-21/4+\varepsilon_3)
     =8(27\varepsilon_3-1)\in\{-8,208\}.
\]
Neither value vanishes for $p\equiv3\pmod4$, $p\geqslant7$.
Hence $x<(P-1)/2$ is impossible for every $p\geqslant7$.

It follows from \eqref{eq:quad-linear-digit} that $v\equiv0\pmod p$.
Since $0\leqslant v<p$, we have $v=0$ and
$T\leqslant(P-1)/2$. Hence
\begin{equation}\label{eq:quad-linear-size}
 \Delta\leqslant\frac{K(P-1)}{K-1}.
\end{equation}
By the bound \eqref{eq:quad-linear-size},
Lemma~\ref{lem:quad-divisibility}~\ref{lift:linear} applies
and gives $p\mid m$. This completes
\eqref{eq:quad-linear-induction}.

\emph{The last step.}
If $f$ is even, the last step has $P=K$, so
$0<\Delta\leqslant K$ and $pK\mid\Delta$, a contradiction.
If $f$ is odd, the last step has $P=pK$. The bound gives
\[
 m\leqslant p+\frac{p-1}{K-1}\leqslant p+1.
\]
Since $p\mid m$, the only positive possibility is $m=p$.
Put $L=pK=\Delta$, so $q=pK^2$ and $r\equiv L+3\pmod h$.
Apply \eqref{eq:quad-block-one} with this $K$ to
\[
 b=\left[-\frac{(K+1)r}{2}\right]_h
  =h-\frac{(p+3)K+2}{2}
  =\frac{pK-p-4}{2}K+\frac{K-1}{2}.
\]
Since $K\geqslant p\geqslant7$, one has
$0\leqslant b\leqslant h-K-2$. Thus $\varepsilon=1$, and
\eqref{eq:quad-linear-digit}, with $K$ in place of $P$, demands
$d+1=-3/2$ in $\mathbb F_p$. But the displayed quotient gives
$d+1=(pK-p-2)/2\equiv-1\pmod p$, a contradiction.
Hence $\Delta=0$.
\end{proof}

\subsection{Split coefficients}\label{subsec:quad-split}

\begin{lemma}\label{lem:quad-split-exponent}
Suppose $p\equiv1\pmod4$, $p\geqslant13$, and
$s^2-2s+5=0$. If $0<j<h/2$ and $F_{r,a}$ permutes
$\mathbb{F}_{q^2}$, then $j=(q-1)/4$.
\end{lemma}
\begin{proof}
Write $s=1+2i$, with $i^2=-1$ in $\mathbb F_p$, and set
\[
 j_0=\frac{q-1}{4},\qquad \Delta=j_0-j\geqslant0,
 \qquad \eta=\frac{\Delta+1/2}{q+1}.
\]
Equation~\eqref{eq:quad-first} gives $p\mid\Delta$, and $0<\eta<1/4$.
If $f=1$, this already gives $\Delta=0$.

We prove the assertion inductively. At a stage $K=p^k\mid\Delta$,
$1\leqslant k\leqslant\lceil f/2\rceil$, suppose either $k=1$
or the preceding stage gave $\eta\leqslant1/[2(K/p-3)]$.
Put $H=(K-1)/4$ and, for $l=1,3$, let
\[
 b_l=H-\left\lceil\frac{l-1}{4}+(K-l)\eta\right\rceil,
 \qquad c=H-b_1,\qquad d=b_1-b_3\in\{0,1\}.
\]
We have $1\leqslant c<p$: for $k=1$,
$c\leqslant H<p$; otherwise, with $L=K/p\geqslant13$,
\[
 (K-1)\eta\leqslant\frac{pL-1}{2(L-3)}
 <\frac{13p}{20}<p-1.
\]

Apply Lemma~\ref{lem:quad-power-sums} with
$m=m_l=(K-l)/2$, $l=1,3$. In neither case is $m_lj$ divisible by $h$:
indeed, $m_lj=b_lh$ would imply $4b_l\equiv l\pmod K$, whereas
$0\leqslant4b_l<K-l$ and $l=1,3$. Since $K\mid\Delta$,
$4j\equiv-1\pmod K$, and the residues modulo $K$ of
$2(m_lj-b_lh)-1$ are exactly
\[
 E_1=2H+c,\qquad E_3=c+d,
 \qquad \sigma_l=(-1)^{E_l}.
\]
Put $\alpha=-s/(1-s)$ and $\beta=s/(1+s)$, which are defined since
$s\ne\pm1$. In \eqref{eq:quad-residue} with $m=m_l$ the prefactor
$(1-s^2)^{2m_l}/(2^{d}(s-1))$ is nonzero, and $d=2m_l-1=K-l-1<K$, so
only the series modulo $z^K$ matters. Since $(1+yz)^K\equiv1\pmod{z^K}$
for every $y\in\mathbb{F}_q$, the factors in \eqref{eq:quad-residue}
satisfy, modulo $z^K$,
\[
\begin{gathered}
 (1+z)^{2(m_lj-b_lh)-1}\equiv(1+z)^{E_l},\qquad
 (1+\alpha z)^{d}\equiv(1+\alpha z)^{-l-1},\\
 (1+\beta z)^{2m_l}\equiv(1+\beta z)^{-l},\qquad
 (1+z/2)^{-2m_l}\equiv(1+z/2)^{l}.
\end{gathered}
\]
Hence the coefficient in \eqref{eq:quad-residue}, which must vanish,
equals
\[
 [z^{K-l-1}]\,\frac{(1+z)^{E_l}(1+z/2)^l}{(1+\alpha z)^{l+1}(1+\beta z)^l}.
\]
Under the substitution $w=s(1+z)$ we have $1-w=(1-s)(1+\alpha z)$,
$1+w=(1+s)(1+\beta z)$, $w^2-s^2=2s^2z(1+z/2)$ and
$(w-s)^K=s^Kz^K$. Therefore this coefficient equals
$C_l\operatorname{Res}_{w=s}\omega_l$, where
$C_l=2^{-l}s^{K-E_l-2l-1}(1-s)^{l+1}(1+s)^l\ne0$ and
\[
 \omega_l=\frac{w^{E_l}(w^2-s^2)^l}
 {(w-s)^K(1-w)^{l+1}(1+w)^l}\,dw.
\]
The residue of $\omega_l$ at infinity is zero since $E_l<K$; at $w=\pm1$,
$(w-s)^{-K}=(w^K-s)^{-1}$ is constant to the required orders
since $s^K=s$ and $l<K$. Summing residues gives
$R_l(E_l,\sigma_l)=0$, where
\begin{align*}
 R_l(X,\sigma)={}&(1+s)[z^l]
 \frac{(1+z)^X((1+z)^2-s^2)^l}{(2+z)^l}\\
 &-\sigma(1-s)[z^{l-1}]
 \frac{(1+z)^X((1+z)^2-s^2)^l}{(2+z)^{l+1}}.
\end{align*}
Only coefficients of orders at most three occur. Since
$E_1\equiv c-1/2$ and $E_3\equiv c+d\pmod p$, their
expansion, using $i^2=-1$, gives
\[
 R_1=\begin{cases}
 8(c-1),&c\text{ odd},\\
 4(2c+i-1),&c\text{ even},
 \end{cases}
\]
and, subject to $R_1=0$, the following values of $R_3$:
\[
\begin{array}{c|cc}
 &d=0&d=1\\ \hline
 c\text{ odd}&0&12(3+i)\\
 c\text{ even}&\frac43(-1+3i)&-\frac83(2+9i).
\end{array}
\]
The norms of $3+i$, $-1+3i$, and $2+9i$ are $10,10,85$.
Thus $R_3=0$ forces $c$ odd and $d=0$, except perhaps when
$p=17$ and $c$ is even. But $i=\pm4$ then gives
$c\in\{7,11\}$ from $2c=1-i$, a contradiction.
Hence $R_1=0$ gives $c=1$, $d=0$, so $b_1=b_3=H-1$ and
\[
 \eta\leqslant\frac1{2(K-3)}.
\]

If $q\leqslant K^2$, this implies
$0\leqslant\Delta\leqslant(q+1)/[2(K-3)]-1/2<K$.
Since $K\mid\Delta$, we are done. Otherwise write $q=KP$
with $P>K$ and $\Delta=mK>0$. The same bound gives
\[
 \Delta\leqslant\frac{KP+1}{2(K-3)}-\frac12
 <\frac{K(P-1)}{K-1},
\]
where the last inequality follows from
$KP(K-5)-K^2+K+4>0$ for $K\geqslant13$ and $P\geqslant K$.
Lemma~\ref{lem:quad-divisibility}~\ref{lift:quadratic} therefore
gives $p\mid m$. This supplies the next induction stage.
After $k=\lceil f/2\rceil$ we have $q\leqslant K^2$,
and the preceding argument gives $\Delta=0$.
\end{proof}

\subsection{Characteristic five}\label{subsec:quad-five}

\begin{lemma}\label{lem:quad-five-exponent}
Suppose $q=5^f$, $s=2$, and $0<j<h/2$. If $F_{r,a}$
permutes $\mathbb{F}_{q^2}$, then
\[
 j=1\quad\text{or}\quad j=\frac{q-1}{4}.
\]
\end{lemma}
\begin{proof}
Equation~\eqref{eq:quad-first} gives $j\equiv1\pmod5$, which proves the
lemma for $q=5$. For $q\geqslant25$ we evaluate the identities of
Lemma~\ref{lem:quad-power-sums} by Lucas' theorem.

For $1\leqslant m<h$, put $d=2m-1$, $R=mj-bh$, and $N=R-1$,
where $b=\lfloor mj/h\rfloor$. In \eqref{eq:quad-binomial-residue},
$s=2$ gives $\lambda=2$ and $s^{2d}=-1$. Hence, writing
$D_d(n)$ for $D_\lambda(d,n)$ with $\lambda=2$,
\[
 C(m,j,b):=-D_d(N)-2D_d(N+h)=\mathbf1_{R=0}.
\]
For $N=-1$ use $N=q-1$, which gives the same coefficients
through degree $d<q$. Writing
$d=\sum d_i5^i$ and $n=\sum n_i5^i$, Lucas' theorem
now gives the product
\begin{equation}\label{eq:quad-five-lucas}
 D_d(n)=\prod_i\Phi_{d_i}(n_i),\qquad
 \Phi_t(v)=\sum_{l=0}^t\binom tl\binom vl2^l,
\end{equation}
where
\[
 \Phi_1(v)=1+2v,\qquad
 \Phi_2(v)=2v^2+2v+1,\qquad
 \Phi_4(v)=(-1)^v\quad(0\leqslant v<5).
\]
Here $\binom4l=(-1)^l$ in $\mathbb F_5$.

\emph{Residues modulo $25$.}
Write $\theta=j/h$ and $\delta=1/2-\theta$. The cases
$m=3,4,11$ of Lemma~\ref{lem:quad-power-sums} give the initial
alternatives
\begin{equation}\label{eq:quad-five-initial}
 \begin{array}{ll}
 j\equiv1\pmod{25},&0<\theta<1/11,\\
 \text{or}\quad j\equiv6\pmod{25},&0<\delta<1/22.
 \end{array}
\end{equation}
Write $j\equiv1+5t\pmod{25}$, $0\leqslant t\leqslant4$.
Formula~\eqref{eq:quad-five-lucas} gives
$C(3,j,b)=2t+3b$ and $C(4,j,b)=b(2t+3)$ for $b=0,1$.
Recall that $C(m,j,b)$ equals $1$ if $mj=bh$ and $0$ otherwise.
For $m=3,4$ this forces
$t=0$ and $\theta<1/4$, or $t=1$ and $\theta>1/3$;
the two equations are incompatible when $4j=h$ or $3j=h$.
Since $d=21=4\cdot5+1$, we have
$D_{21}(n)=(1+2n_0)(-1)^{n_1}$ for the two low
base-five digits of $n$. In the surviving intervals,
\[
 C(11,j,b)=b(4-b)\ (b=0,1,2),\qquad
 C(11,j,b)=2(5-b)\ (b=3,4,5).
\]
Since $C(11,j,b)$ equals $1$ if $11j=bh$ and $0$ otherwise,
this forces $b=0$, or $b=5$ with $11j>5h$.
This proves \eqref{eq:quad-five-initial}.

\emph{Induction.}
Let $K=5^k$, $2\leqslant k\leqslant f$. We prove inductively that
\begin{equation}\label{eq:quad-five-induction}
 \begin{array}{ll}
 \mathrm L:\ j\equiv1\pmod K,&0<K\theta\leqslant5/2,\\
 \mathrm H:\ j\equiv(K-1)/4\pmod K,&0<K\delta\leqslant3/2.
 \end{array}
\end{equation}
The base case $K=25$ follows from \eqref{eq:quad-five-initial}.
Suppose $5K\leqslant q$ and use $m=2K+3$; then
$m<h$ and $d=4K+5<q$. In case L put
$j\equiv1+uK\pmod{5K}$, $b=n$, with
$0\leqslant u\leqslant4$ and $0\leqslant n\leqslant5$.
In case H put
$j\equiv(K-1)/4+uK\pmod{5K}$,
$b=K+1-n$, with $0\leqslant n\leqslant3$.
In neither case is $mj=bh$: reducing $2mj=b(q+1)$
modulo $K$ would give $n\equiv6$ in L or $2n\equiv5$
in H, contrary to the stated integer ranges.

For $n\geqslant1$, put $\varepsilon=(-1)^{[2-2u]_5}$,
$\varepsilon'=(-1)^{[4-2u]_5}$, with residues in
$\{0,1,2,3,4\}$, and $\nu=[n]_2$. The only nontrivial
base-five digits of $d=4K+5$ are at positions $1,k$.
Using $2N\equiv2mj-b(q+1)-2\pmod{5K}$, the values
$C(m,j,b)$ in cases L and H, respectively, are
\begin{equation*}
 C_{\mathrm L}=-(1+\nu)\varepsilon,\qquad
 C_{\mathrm H}=(2-3\nu)\varepsilon+(1-3\nu)\varepsilon'.
\end{equation*}
Neither value can vanish: they are a nonzero signed multiple
of $1$ or $2$, or a signed sum of $1$ and $2$ in $\mathbb F_5$.
The required vanishing therefore forces $n=0$ in both cases, and
hence
\[
 0<\theta<\frac1{2K+3}\quad(\mathrm L),\qquad
 0<\delta<\frac1{2(2K+3)}\quad(\mathrm H).
\]

To determine $u$, use $m=K+1$, $b=0$ in L and
$m=K-1$, $b=(K-3)/2$ in H. In L, $d=2K+1$
has digits $1,2$ at positions $0,k$; the congruences
$N\equiv(1+u)K$ and $h\equiv(5K+1)/2\pmod{5K}$
give, by \eqref{eq:quad-five-lucas},
\[
 C(K+1,j,0)=\Phi_2(u+3)-\Phi_2(u+1)=3u.
\]
In H, $d=2K-3$ has digits $2,4,\ldots,4,1$;
$N\equiv(3-u)K\pmod{5K}$, while the units digit
of $N+h$ is $3$ and $\Phi_2(3)=0$. Thus
\[
 C\bigl(K-1,j,(K-3)/2\bigr)=3+2u.
\]
The strict bounds above show that $b=\lfloor mj/h\rfloor$ and
$mj\ne bh$ for these choices, so both values must vanish; this forces
$u=0$ in L and $u=1$ in H. Hence L lifts to
$j\equiv1\pmod{5K}$ and H to $j\equiv(5K-1)/4\pmod{5K}$,
with \eqref{eq:quad-five-induction} at $5K$.
At $K=q$, since $0<j<h/2<q$, the two congruences are exactly
$j=1$ and $j=(q-1)/4$.
\end{proof}

\subsection{Proof of the quadratic classification}\label{subsec:quad-proof}

\begin{proof}[Proof of Theorem~\ref{thm:quadratic}]
\emph{Necessity.} Suppose that $F_{r,a}$ permutes $\mathbb{F}_{q^2}$.
Lemma~\ref{lem:coefficient-product} gives $\gcd(r,q-1)=1$ and
$s\in\mathbb{F}_q^\times$. If $s^2=1$,
Lemma~\ref{lem:quad-elementary-families}~\ref{quad-family:norm-one} gives
condition~(i). Assume $s^2\ne1$. If $j=0$,
Lemma~\ref{lem:quad-elementary-families}~\ref{quad-family:zero} gives
$s^2=-1$, the second row of Table~\ref{tab:quadratic}. Otherwise,
by \eqref{eq:quad-reciprocity} we may assume $0<j\leqslant h/2$, since
this reciprocity interchanges the first and third rows and the fourth
and fifth rows, and maps the second and sixth rows to themselves. If
$j=h/2$, Proposition~\ref{prop:quad-coefficient} gives $q=11$ and
$s\in\{2,6\}$, the sixth row. If $0<j<h/2$, the same proposition leaves
$3s=1$ with $p\ne3$, or $s^2-2s+5=0$:
\begin{itemize}[leftmargin=*]
\item for $p=3$, only $s^2-2s+5=0$ is possible, and
Proposition~\ref{prop:quad-nonsplit} gives the fourth row;
\item for $p=5$, both conditions give $s=2$, which satisfies $3s=1$ and
$s^2-2s+5=0$, and Lemma~\ref{lem:quad-five-exponent} gives the third or
the fourth row;
\item for $p\geqslant7$ and $3s=1$, Lemma~\ref{lem:quad-linear-exponent}
gives $j=1$, the third row;
\item for $p\geqslant7$ and $s^2-2s+5=0$,
Proposition~\ref{prop:quad-nonsplit} if $p\equiv3\pmod4$, or
Lemma~\ref{lem:quad-split-exponent} if $p\equiv1\pmod4$ (so that
$p\geqslant13$), gives $j=(q-1)/4=(h-1)/2$ and $q\equiv1\pmod4$, the
fourth row.
\end{itemize}

\emph{Sufficiency.} Assume $\gcd(r,q-1)=1$. Condition~(i) is
Lemma~\ref{lem:quad-elementary-families}~\ref{quad-family:norm-one}; the
first three rows are
Lemma~\ref{lem:quad-elementary-families}~\ref{quad-family:zero}
and~\ref{quad-family:linear}; the fourth and fifth rows are
Proposition~\ref{prop:quad-middle}; and the sixth row is the converse
part of Proposition~\ref{prop:quad-coefficient}.
\end{proof}

\section{Extensions of even degree at least four}\label{sec:higher}

In this section we prove Theorem~\ref{thm:higher}; the argument is
completed in Section~\ref{subsec:higher-proof}. Throughout, $q$ is odd
and $D:=2(q-1)$. The coefficient restriction in
Lemma~\ref{lem:coefficient-product} allows us to descend whenever the
extension degree is divisible by four (Section~\ref{sec:higher-descent}).
The degrees $e=2u$ with $u\geqslant3$ odd are treated in
Sections~\ref{sec:higher-exponent}--\ref{sec:higher-ternary}.

\subsection{Descent and degree four}\label{sec:higher-descent}

\begin{lemma}\label{lem:higher-descent}
Suppose that $4\mid e$ and that $X^r(X^D+a)$ permutes
$\mathbb{F}_{q^e}$, where $a\ne0$. Then there exists
$b\in\mathbb{F}_{q^{e/2}}^\times$ such that $X^r(X^D+b)$
permutes $\mathbb{F}_{q^{e/2}}$.
\end{lemma}
\begin{proof}
Put $Q=q^{e/2}$. As $Q\equiv1\pmod D$,
$\gcd(Q+1,D)=2$. The cyclic group $\mathbb F_{Q^2}^\times$
therefore satisfies
\[
 \mathbb F_Q^\times(\mathbb F_{Q^2}^\times)^D
 =(\mathbb F_{Q^2}^\times)^2.
\]
By Lemma~\ref{lem:coefficient-product}, $a$ is a square;
write $a=bc^D$ with $b\in\mathbb F_Q^\times$ and
$c\in\mathbb F_{Q^2}^\times$. Then
\[
 c^{-r-D}F_{r,a}(cX)=X^r(X^D+b)
\]
permutes $\mathbb F_{Q^2}$. Its restriction to $\mathbb F_Q$
is injective and hence bijective.
\end{proof}

\begin{proposition}\label{prop:higher-quartic}
For every $r\geqslant1$ and every $a\in\mathbb{F}_{q^4}^\times$,
the binomial $X^r(X^D+a)$ does not permute $\mathbb{F}_{q^4}$.
\end{proposition}
\begin{proof}
Suppose $F_{r,a}$ permutes $\mathbb F_{q^4}$. Put
$c=(q+1)/2$, $N=q^2+1$, $\lambda=cN=(q^4-1)/D$,
and $w=a^\lambda$. Lemma~\ref{lem:coefficient-product} gives
$w\in\mathbb F_q^\times$; since $\lambda$ is even, $w\ne1$.
Choose $b\in\mathbb F_{q^2}^\times$ with $b^{q+1}=w$.
Then $b^\lambda=w$, so scaling and restriction yield
permutations $X^r(X^D+b)$ and $X^r(X^D+b\zeta)$ of
$\mathbb F_{q^2}$, where $\zeta$ generates $\mu_{q+1}$
(and $\zeta^\lambda=1$). As $\zeta^c=-1$, their
parameters $s=(-b)^c$ have opposite signs.
Lemma~\ref{lem:quad-first} forces $r\equiv2\pmod c$
and $w=s^2=-1$. For $q=3$, $r$ would be even.

Let $q\geqslant5$. In the power sum of exponent
$AN=q^4-1-DN$, where $A=(q-1)^2=(q-2)q+1$, a
contributing binomial index $i$ satisfies $2i\equiv0\pmod N$.
If $i$ is an odd multiple of $N/2$, its two base-$q^2$
blocks have units digits differing by $\pm c\pmod q$; Lucas' theorem
excludes this, since both base-$q^2$ blocks of $AN$ have units digit $1$.
Thus $i=Nt$. Lucas' theorem gives $\binom{AN}{Nt}=\binom At^2$;
writing $t=t_0+qt_1$, nonvanishing requires
$t_0\in\{0,1\}$ and $0\leqslant t_1\leqslant q-2$.
Since $t\equiv r\equiv2\pmod c$ and $q\equiv-1\pmod c$,
only $t=q(c-2)$ and $t=1+q(c-1)$ remain.
They sum to $A$, and their common coefficient is nonzero:
$\binom{q-2}{c-2}=(-1)^{c-2}(c-1)\ne0$ in $\mathbb F_p$.
The powers of $a$ differ by $2\lambda=(q+1)N$;
vanishing would require $1+a^{-2\lambda}=0$, contrary
to $a^\lambda=-1$.
\end{proof}

\subsection{An exponent congruence in degree \texorpdfstring{$2u$}{2u}}\label{sec:higher-exponent}

Throughout this subsection $q=p^f$ is an odd prime power, $u\geqslant3$ is
odd, and $\ell_i=(q^i-1)/(q-1)$ as in Section~\ref{sec:intro}. Put
\[
 Q:=q^u,\quad P:=q^{u-1},\quad d:=q-1,
 \quad S:=\ell_{u-1},\quad E:=\ell_u=P+S=qS+1,
\]
$K:=\mathbb{F}_Q$ and $L:=\mathbb{F}_{Q^2}$; assume $a\in L^\times$, and
write $F:=X^r(X^{2d}+a)$. By \eqref{eq:reciprocity} with $M=Q^2$ and
$D=2d$, replacing $(r,a)$ by $(r',a^{-1})$, where
$r'\equiv-r-2d\pmod{Q^2-1}$ and $r'\geqslant1$, preserves the
permutation property; this is the inversion used below.

\begin{lemma}\label{lem:higher-digit-word}
Let $q=p^f$ be odd, $u\geqslant3$ odd, and
$r\equiv-2qj\pmod E$, where $1\leqslant j\leqslant T=(P-1)/(q^2-1)$ and
$j\equiv1\pmod q$. For $1\leqslant s\leqslant u-2$, write
$A_s=1+\ell_{s+1}$ and $A_sj=n_s E+R_s$, $0\leqslant R_s<E$.
Suppose that for each such $s$ we have
\begin{equation}\label{eq:higher-digit-conditions}
  1\leqslant R_s\leqslant2S+q^s,\qquad
  \left\lfloor\frac{(n_s+1)E-A_sj}{q^s}\right\rfloor
                             \equiv-1\pmod p ,
\end{equation}
and that $[qj\ell_g]_E\geqslant q^g$ for every $1\leqslant g\leqslant u-2$.
Then $r\ell_t\equiv2\pmod E$ for some
$1\leqslant t<u$ with $\gcd(t,u)=1$.
\end{lemma}
\begin{proof}
\emph{The digits of $j$.}
Write $j=\sum_{b=0}^{u-1}J_bq^b$. Then $J_{u-1}=J_{u-2}=0$
and $J_0=1$. Set $M=(q-1)E=q^u-1$, $\delta=(q-2)j<P/q$,
and $Y_s=[q^{s+1}j]_M$. Since $(q-1)A_s=q^{s+1}+q-2$,
\begin{equation*}
 [Y_s+\delta]_M=(q-1)R_s
 \leqslant 2(P-1)+(q-1)q^s<3P.
\end{equation*}
If $J_b\geqslant2$ is the first nonbinary digit, choose $s=u-2-b$.
Then $Y_s\geqslant2P+q^{s+1}>2(P-1)+(q-1)q^s$, so
$Y_s+\delta$ must wrap modulo $M$. But $\delta<P/q$ makes
both leading digits of $Y_s$ equal to $q-1$, forcing
$J_{b-1}=q-1$, a contradiction. Thus every $J_b\in\{0,1\}$;
there is no wrap, and $n_s=\lfloor j/q^{u-s-1}\rfloor$.

\emph{Symmetry of the digits.}
To prove $J_s=J_{u-1-s}$, induct on $1\leqslant s\leqslant(u-3)/2$.
Assuming the previous equalities, let $c=J_{u-1-s}$ and $x=J_s$.
The coefficient of $q^s$ in $(n_s+1)E-A_sj$ is $c-2x$;
its lower part is $c\ell_s-1-\sum_{l=1}^{s-1}J_lw_l$, where
$w_l=((q-2)q^l+q^{s-l})/(q-1)>0$ and
$\sum_l w_l=\ell_s-1$. Thus the lower part lies in
$[-\ell_s,-1]$ for $c=0$ and $[0,\ell_s)$ for $c=1$.
The carry is $c-1$, and the second condition in \eqref{eq:higher-digit-conditions}
reads $2c-2x-1\equiv-1\pmod p$. Thus $c=x$ and
$W=0v1$ has a palindromic interior $v$.
Write $W'=1v0$, of value $j+P-1$.

\emph{Excluded rotations.}
For a left rotation $Y$ of $W$ by $t$, $2\leqslant t<u$, write
$y$ for its base-$q$ value and $z=t-1$. The last hypothesis
forbids $Y\leqslant W$ with the last $z$ digits of $W$, or $Y>W'$
with the last $z$ digits of $W'$. Indeed, in either case
\[
 i=\left[\frac{j-y}{q-1}\right]_E\in q^z\mathbb Z\cap[0,P),
 \qquad i\equiv-\ell_tj\pmod E.
\]
Rotation preserves digit sums; suffix agreement gives $q^z\mid i$,
using $M\equiv-1$ and $W'\equiv W-1\pmod{q^z}$ in the second case.
As $q^{z+1}\ell_{u-z-1}=E-\ell_t$, we get
$[qj\ell_{u-z-1}]_E=i/q^z<q^{u-z-1}$, a contradiction.

\emph{Comparison with lines of irrational slope.}
For irrational $0<\alpha<1$, let
$v(\alpha)_i=\lfloor(i+1)\alpha\rfloor-\lfloor i\alpha\rfloor$
($1\leqslant i\leqslant u-2$). Its prefix sums $\lfloor k\alpha\rfloor$
($2\leqslant k<u$) make $v(\alpha)$ lexicographically increasing,
with endpoint values $0^{u-2}$ and $1^{u-2}$.
It changes only at reduced rationals $a/t$, $2\leqslant t<u$.
If $v$ is missed, it is skipped at such a breakpoint, whose
periodic limits are $(c01)^\infty$ and $(c10)^\infty$:
only multiples of $t$ change there. Here $|c|=t-2$ and $c$ is palindromic:
$\lfloor a-x\rfloor=a-\lceil x\rceil$ identifies its
$i$-th and $(t-1-i)$-th digits.
Write $v=cbR$. Palindromicity makes $c$ a suffix of $v$ too, so $W=0v1$
ends with $c1$ and $W'=1v0$ ends with $c0$. For a finite nonempty word
$X$ let $X^\infty=XXX\cdots$ denote its periodic extension; finite words
of equal length, and infinite words, are compared lexicographically.
For finite nonempty words $X$ and $Y$,
\[
 X^\infty<Y^\infty\iff XY<YX.
\]
Indeed, if $XY<YX$, then interchanging adjacent factors $YX$ one at a
time gives $X^{|Y|}Y^{|X|}<Y^{|X|}X^{|Y|}$; as $X^{|Y|}$ and $Y^{|X|}$
have the same length, this forces $X^{|Y|}<Y^{|X|}$, whence
$X^\infty=(X^{|Y|})^\infty<(Y^{|X|})^\infty=Y^\infty$. If $XY=YX$, the
same interchanges give $X^{|Y|}=Y^{|X|}$ and $X^\infty=Y^\infty$; the
case $XY>YX$ is symmetric. Taking $X=AZ$ and $Y=A$, and cancelling the
common prefix $A$ of $(AZ)A$ and $A(AZ)$, we obtain
\[
 (AZ)^\infty<A^\infty\iff ZA<AZ,
\]
and likewise with both inequalities reversed; here $ZA$ is the left
rotation of $AZ$ by $|A|$.

If $b=1$, then $W=AZ$ with $A=0c1$ and $Z=R1$, and
$W^\infty<(0c1)^\infty$. Hence the $t$-rotation $ZA=R10c1$ of $W$
satisfies $ZA<W$ and has the same last $t-1$ digits $c1$ as $W$. If
$b=0$, then $W'=A'Z'$ with $A'=1c0$ and $Z'=R0$, and
$(W')^\infty>(1c0)^\infty$, so $R01c0=Z'A'>W'$. The $t$-rotation of
$W=(0c0)(R1)$ is $R10c0>R01c0>W'$, and its last $t-1$ digits $c0$ are
those of $W'$. Both cases contradict the rotation exclusion.
Thus $v=v(\alpha)$ for some irrational $\alpha$.

\emph{Conclusion.}
Let $f_k=\lfloor k\alpha\rfloor$, $h=1+f_{u-1}$.
The palindrome gives $f_k+f_{u-k}=h-1$, and the defining
floor inequalities yield
\[
 f_k<kh/u<f_k+1\qquad(1\leqslant k<u).
\]
Consequently $\gcd(h,u)=1$, and $W$, $W'$ are the lower and
upper words of slope $h/u$. For $1\leqslant t<u$ with
$th\equiv-1\pmod u$, the identity
$\lfloor(m-1)/u\rfloor=\lceil m/u\rceil-1$ for integers $m$
shows that the left $t$-rotation sends $W$ to $W'$.
Hence $(q^t-1)j\equiv P-1\pmod M$, or
$qj\ell_t\equiv-1\pmod E$. Since $r\equiv-2qj\pmod E$,
we obtain $r\ell_t\equiv2\pmod E$.
\end{proof}

\begin{proposition}\label{prop:higher-exponent}
Suppose that $F$ permutes $L$. After inversion if necessary, there is
an integer $t$ with $1\leqslant t<u$ and $\gcd(t,u)=1$ such that
\[
 r\ell_t\equiv2\pmod E.
\]
\end{proposition}

\begin{proof}
\emph{Normalization.}
Lemma~\ref{lem:coefficient-product} makes $a$ a square. If $q=3$,
we may further take $a=\gamma$, an element of order four in
$\mathbb{F}_9$. Indeed, a permutation would give
$a^{(Q^2-1)/4}=\pm1$, and the positive sign would make $-a$
a fourth power and hence give a nonzero root. Since $Q^2\equiv9\pmod{16}$, an element $\gamma\in\mathbb{F}_9$
of order four has the same non-fourth-power class; scaling
yields $a=\gamma$, with $\gamma^Q=\gamma^{-1}$.

\emph{The basic identity.}
For $1\leqslant B<E$ and $m=(q-1)B<(Q+1)/2$, put
$I(B)=\{0\leqslant i\leqslant m:2i\equiv-rB\pmod E\}$ and $w=a^{Q+1}$.
The vanishing power sum of $F^{m(Q+1)}$ gives
\begin{equation*}
 \sum_{i\in I(B)}\binom mi^2w^{m-i}=0.
\end{equation*}
Indeed, the two binomial indices $i,l\leqslant m$ satisfy
$(Q+1)/2\mid i-l$, so $i=l$. Such a sum cannot have just one
or two nonzero terms. For $q\equiv3\pmod4$, each term is a
square in $K$ and $-1$ is not; for $q=3$ these terms equal
$1$, since $w=\gamma^{Q+1}=1$. For $q\equiv1\pmod4$,
scaling takes the square $a$ into $K^\times$ because
$\gcd(Q+1,2(q-1))=2$; restriction then gives both
$\sum_{I(B)}X_i=0$ and $\sum_{I(B)}X_i^2=0$, where
$X_i=\binom mi a^{m-i}$. Two nonzero terms cannot satisfy
both equalities in odd characteristic.

\emph{An additional identity for $q=3$.}
For $q=3$, take also $B=2S-\ell_s$, $m=2P-1-3^s$
($0\leqslant s\leqslant u-2$) and $H=(Q+1)/2=E+1$.
The full power sum has terms $\binom mi\binom ml\gamma^{l-i}$ with
$E(Q+1)\mid rB(Q+1)+2(i+Ql)$. Since $m<2H$,
these are diagonal ($i=l$, $2i\equiv-rB\pmod E$) or
pairs $(i,l)=(b+H,b),(b,b+H)$ with
$2b+1\equiv-rB\pmod E$. The bounds
$m-E=S-3^s$ and $m-H=S-1-3^s$ allow at most two
diagonals and one pair. Two diagonals cannot coexist with a pair:
the smaller diagonal $i$ and a pair index $b$ would give
$2(i-b)\equiv1\pmod E$, yet $|2(i-b)|\leqslant2(S-3^s)<E-1$.
Lucas' theorem gives coefficient $1$ for each surviving
diagonal, since the ternary digits of $m$ are $1$ at $s,u-1$
and $2$ elsewhere. For a surviving pair, the two binomial
coefficients given by Lucas' theorem multiply to $1$ (the top digits are $0,1$,
and the $s$-digits sum to $1$); as $H\equiv2\pmod4$,
the pair contributes $-2=1$ in $\mathbb F_3$.
Thus the nonempty sum has at most two contributions, each $1$,
contradicting its vanishing. In particular all diagonal terms vanish.

\emph{The residue of $r$ modulo $E$.}
At $B=S$ we have $m=P-1<E$, with every binomial coefficient
nonzero. The missing residue class must be $E-j$ with
$1\leqslant j\leqslant S$. Since $qS=E-1$,
\begin{equation*}
 r\equiv-2qj\pmod E,
 \qquad i\equiv qjB\pmod E\quad(i\in I(B)).
\end{equation*}
We need $j\equiv1\pmod q$. For $q=3$, the case $s=0$ of the preceding
paragraph forces the diagonal index $E-2j$ to have its last ternary
digit equal to $2$, whence $j\equiv1\pmod3$.
For $q\geqslant5$, suppose $j\not\equiv1\pmod q$ and let
$0\leqslant k<f$ be the first nonzero base-$p$ digit of $j-1$. Set
$t=p^k$, $B=(t+1)S$, $m=(t+1)(P-1)$,
$A=(t+1)j<P$. The possible indices are
$i_l=lE-A$ ($1\leqslant l\leqslant t$), with $m<(Q+1)/2$.
As $m=tP+(P-1-t)$, Lucas' theorem permits only
$i_1$ if $A>S$ or $i_t$ if $A\leqslant tS$; at least one
is in range. The residue $P-1-t$ of $m$ modulo $P$ has base-$p$ digit
$p-2$ at position $k$, so an index whose digit there is $p-1$ has zero
coefficient. For $k=0$, $i_1=E-2j$ survives:
its units digit $1-2j$ differs from $-1\pmod p$.
For $k\geqslant1$ write $j=1+ct$, $p\nmid c$; then
\[
 i_1\equiv-(c+1)t,\qquad i_t\equiv-1-ct\pmod{pt}.
\]
Both endpoints have $p^k$-digit $-c-1\ne-1\pmod p$,
so every in-range endpoint survives. The resulting one
or two nonzero terms contradict the fact, shown above, that the sum
cannot have exactly one or two nonzero terms.
Hence $j\equiv1\pmod q$.

\emph{Reduction to small $j$.}
Inversion replaces $j$ by $S+1-j$. At $B=S-1$,
$m=P-q$; any $S<(q+1)j\leqslant E$ would give the unique
nonzero index $E-(q+1)j$. The two remaining intervals
are interchanged by inversion. Hence, writing
$T=(P-1)/(q^2-1)$, we may assume
\begin{equation*}
 1\leqslant j\leqslant T,\qquad j\equiv1\pmod q.
\end{equation*}
For $u=3$ we have $T=1$, so inversion gives
$r\equiv2\pmod E$.

\emph{Digit conditions for $u\geqslant5$.}
For $u\geqslant5$, take $B=2S-\ell_s$, $m=2P-1-q^s$ and
$A_s=1+\ell_{s+1}$ ($1\leqslant s\leqslant u-2$).
There are at most two indices $i\equiv-A_sj\pmod E$;
the residue $P-1-q^s$ of $m$ modulo $P$ has base-$p$ digit $p-2$ at
position $fs$, so an index whose digit there is $p-1$ has zero
coefficient.
Indices differing by $E$ have digits there differing by
$1$ or $2\pmod p$, so both coefficients cannot vanish.
The arguments above excluding sums with exactly one or two nonzero
terms (for $q=3$, the paragraph on $B=2S-\ell_s$) therefore force the
sole index $E-R_s$ to have zero coefficient. Its range and missing digit give
\eqref{eq:higher-digit-conditions}, with $R_s\ne0$.
For $B=S-\ell_z$, $1\leqslant z\leqslant u-2$, we have $m=P-q^z<E$ and
$i\equiv-\ell_{z+1}j\pmod E$. A multiple of $q^z$ in the
allowed range would contribute a single nonzero term, by Lucas'
theorem.
Since $q^{z+1}\ell_{u-z-1}=E-\ell_{z+1}$, its absence yields
$[qj\ell_{u-z-1}]_E\geqslant q^{u-z-1}$.
Apply Lemma~\ref{lem:higher-digit-word}.
\end{proof}

\subsection{Norm ratios on the projective line}\label{sec:higher-norm}

This subsection supplies the counting input for the collision argument
of Section~\ref{sec:higher-collision}.
Write $k:=\mathbb{F}_{q^2}$, $N:=N_{k/\mathbb{F}_q}$,
$U:=\ker N$ and $G:=k^\times/U^2$.
For $b,b^q,d,d^q\in k$ pairwise distinct and $n$ odd define
\[
 W_n(\lambda):=N_{\mathbb{F}_{q^{2n}}/k}
    \!\left(\frac{\lambda-b}{\lambda-d}\right)
 \quad(\lambda\in\mathbb{F}_{q^n}),\qquad W_n(\infty):=1 .
\]
The numerator and denominator are nonzero because $n$ is odd.
Brackets denote classes in $G$.

\begin{proposition}\label{prop:higher-norm}
Let $q\geqslant5$ be odd, and let $b,b^q,d,d^q$ be pairwise distinct.
For every odd $n\geqslant5$ and every $g\in G$, at least $n$ points
$\lambda\in\mathbf P^1(\mathbb{F}_{q^n})$ satisfy $[W_n(\lambda)]=g$.
If, in addition, $N(b)=N(d)$ is a nonsquare in $\mathbb{F}_q$
and $d/b\in U^2$, then $[W_3]$ is surjective onto $G$;
every class $g$ of nonsquare norm occurs at least $q+4$ times.
\end{proposition}

\begin{proof}
Put $S=\mathbf P^1(\mathbb{F}_q)$ and write $[z]$ for the class of
$z\in k^\times$ in $G$. For $x\in S$ let
$w_x=(x-b)/(x-d)$, with $w_\infty=1$; these $q+1$ weights are distinct.
For a homogeneous form $H$ over $\mathbb{F}_q$, up to nonzero scalars,
set $w(H)=H(b,1)/H(d,1)$ whenever both evaluations are nonzero.

\emph{Weights of quadratic forms.}
The quadratic equation $H(b,1)=sH(d,1)$, $s\in k^\times$,
has a nonzero solution: it imposes two linear conditions on three
coefficients. Its kernel has dimension one or two, and no nonzero
member vanishes at $b$ or $d$, since their irreducible quadratics
are distinct. Of the $q^2+q-1$ admissible projective quadratics,
exactly one weight $h$ occurs $q+1$ times; every other weight occurs
once. The exceptional pencil has no common zero in $S$, since
otherwise all linear forms would have weight $h$. Its roots induce
an involution $x\mapsto x^*$ of $S$ with $w_xw_{x^*}=h$.

\emph{A weighted Euler product.}
In $R=\mathbb Q[G]$ put $v_x=[w_x]$, $A=\sum_{x\in S}v_x$,
$C=[h]$, $J=\sum_{g\in G}g$, and $B(T)=1-AT+qCT^2$.
Unique factorization gives a weighted Euler product for projective
forms,
\[
 Z(T)=\sum_H[w(H)]T^{\deg H}
 =\prod_P(1-[w(P)]T^{\deg P})^{-1}
 \equiv B(-T)\pmod{\mathbb QJ}.
\]
Indeed, the quadratic count gives $[T^2]Z=qC\pmod{\mathbb QJ}$;
for degrees $m\geqslant3$, evaluation at the two distinct irreducible
quadratics and the Chinese remainder theorem make the weights
uniform. The product includes the irreducible form $Y$.

\emph{Cyclic words.}
Call a cyclic word $x_1\cdots x_n$ on $S$ \emph{admissible} if
$x_{i+1}\ne x_i^*$, including the cyclic edge. Set
$D_{xy}=v_y\mathbf1_{\{y=x^*\}}$ and $M=\mathbf1v-D$, where
$v=(v_x)_{x\in S}$. Since $D^2=CI$ and
$vD\mathbf1=(q+1)C$, the matrix determinant lemma gives
\[
 \det(I-TM)=\det(I+TD)\frac{B(T)}{1-CT^2}.
\]
Now $\operatorname{tr}(M^n)$ counts admissible cyclic words with
chosen initial position, and
$\operatorname{tr}(D^n)=\sum_{x=x^*}v_x^n$ for odd $n$.
Logarithmic derivatives therefore yield
\[
 \operatorname{tr}(M^n)+\sum_{x=x^*}v_x^n
 =-[T^n]\frac{TB'(T)}{B(T)}
 =[T^n]\frac{TZ'(T)}{Z(T)} \quad(n\text{ odd}),
\]
initially modulo $\mathbb QJ$. Both the cyclic-word and
Euler-product coefficients have augmentation $q^n+1$, so equality
holds in $R$. The Euler product counts each degree-$j$ point
($j\mid n$) with weight $[w(P)]^{n/j}$, $j$ being odd. Hence
\begin{equation}\label{eq:higher-point-word}
 \sum_{\lambda\in\mathbf P^1(\mathbb{F}_{q^n})}[W_n(\lambda)]
  =\operatorname{tr}(M^n)+\sum_{x=x^*}v_x^n
 \qquad(n\text{ odd}).
\end{equation}
In particular, every admissible cyclic word gives a class
attained by $W_n$.

\emph{The letter weights.}
Put $c=(b-d^q)/(d-d^q)$ and
$\delta=(b-b^q)/(d-d^q)\in\mathbb F_q^\times$.
The condition $x^q=x$ on $w_x=(x-b)/(x-d)$ gives
$\{w_x:x\in S\}=cZ_0$, where
\[
 Z_0=\{z\in k:N(z-1)=t\},\qquad
 t=1-\delta/N(c)\notin\{0,1\}.
\]
For $s=N(z)$, the elements with norm $s$ are roots of
$X^2-(s+1-t)X+s$. Thus $V=N(Z_0)$ has at least
$(q+1)/2$ elements and contains a nonsquare $s_0$.
Its conjugate lifts $z_0,z_0^q$ represent both classes
above $s_0$ in $G$, since $z_0^{q-1}\notin U^2$.

\emph{Odd $n\geqslant5$.}
Fix odd $n=3+2k\geqslant5$ and $g\in G$, and put
$g'=gC^{-k}$. Since $2|V|>q-1$, choose $s_1,s_2\in V$ with
$s_1s_2=N(g'/[c]^3)/s_0$ and $(s_1,s_2)\ne(s_0,s_0)$.
Lift them to $z_1,z_2\in Z_0$, and choose $z_3$ among
$z_0,z_0^q$ to give $[c^3z_1z_2z_3]=g'$.
The three letters are not all equal, hence can be ordered
$x,y,z$ with $x\ne z^*$. Inserting $a,a^*$ around $y$,
with $a\notin\{x^*,y^*,y,z\}$, gives an admissible five-word
of class $g'C$; such $a$ exists since $q+1>4$.
Repeating $k-1$ times gives an admissible $n$-word of class
$g$. Since the inserted letters avoid $y$ and $y^*$, the
original $y$ occurs at most twice; oddness of $n$ makes the
word primitive. Its $n$ rotations give the claimed lower
bound by \eqref{eq:higher-point-word}.

\emph{The case $n=3$.}
For $n=3$ assume the extra hypotheses, and put $A_0=1-t$.
They imply
\[
 t=N(b)\left(\frac{b-d}{bd-N(b)}\right)^2
 \notin(\mathbb F_q^\times)^2,\qquad c^{q-1}=d/b\in U^2.
\]
The fraction lies in $\mathbb F_q$, so $N(c)$ is square.
Now $z+z^q=N(z)+A_0$ implies $N(z-A_0)=tN(z)$:
$z\mapsto A_0/z$ is an involution of $Z_0$ with at most
two fixed points. Since $|Z_0|=q+1\geqslant6$, two distinct
nonfixed orbits yield quadratic forms of common weight
$c^2A_0$. Thus $h=c^2A_0$, and $N(C)$ is square.

For a class $g$ of square norm, choose $s_1,s_2\in V$
with $s_1s_2=N(g/[c]^3)/s_0$ and lift to a triple of
weight $g$. If it is admissible, we are done. Otherwise a
forbidden adjacent pair has product $A_0$; the remaining
letter $z$ has square norm. Hence $[z]=[z^q]$, although
$z\ne z^q$ because $t$ is nonsquare. The triple
$z,z,A_0/z^q$ has the same weight and is admissible
unless $z^2=A_0$. In that case $z$ is a fixed letter
$x\in S$, and $[W_3(x)]=[cz]^3=g$.

For a class $g$ of nonsquare norm, every $s\in V$ has
two conjugate lifts, so $|V|=(q+1)/2$; above any nonsquare
$s$, the lifts represent distinct classes of $G$. For each
$s_3\in V$, at least two pairs $s_1,s_2\in V$ satisfy
$s_1s_2s_3=N(g)/N(c)^3$. Each norm triple has eight
lifts, four of class $g$, since conjugating a letter of
nonsquare norm exchanges the two classes above the product.
Thus $[A^3]_g\geqslant4(q+1)$; as $C$ has square norm,
$[AC]_g\leqslant1$. The cubic case of \eqref{eq:higher-point-word}
gives at least $4(q+1)-3q=q+4$ points of class $g$.
\end{proof}

\begin{corollary}\label{cor:higher-norm-three}
Let $q=3$, let $b$ generate $\mathbb{F}_9^\times$, and set $d=-b$. For every
odd $n\geqslant3$, some $\lambda\in\mathbf P^1(\mathbb{F}_{3^n})$ has $[W_n(\lambda)]$
equal to the unique element of order two in $\mathbb{F}_9^\times/\{\pm1\}$.
\end{corollary}
\begin{proof}
The letter weights are $1,-1,w,w^{-1}$, where
$w=(1-b)/(1+b)$ has norm $-1$ and order eight.
Since $X^2,Y^2$ both have weight $1$, the exceptional weight
is $h=1$, and partners have inverse weights. For odd $n\geqslant3$,
the cyclic word $w,w,1,-1,1,\ldots,1$ (with $n-2$ alternating
entries after the two $w$'s) is admissible and has class $[w]^2$,
the unique class of order two in $\mathbb F_9^\times/\{\pm1\}$.
Apply \eqref{eq:higher-point-word}.
\end{proof}

\subsection{Collisions}\label{sec:higher-collision}

\begin{proposition}\label{prop:higher-odd}
Let $q\geqslant5$ be an odd prime power and $u\geqslant3$ an odd integer.
For every $r\geqslant1$ and $a\in\mathbb{F}_{q^{2u}}^\times$,
$X^r(X^{2(q-1)}+a)$ does not permute $\mathbb{F}_{q^{2u}}$.
\end{proposition}

For any odd $q$ and odd $u\geqslant3$, we use $Q$, $K$, $L$, $E$ and $F$
as in Section~\ref{sec:higher-exponent}, $k$ and
$U=\ker N=\mu_{q+1}$ as in Section~\ref{sec:higher-norm}, and
$D=2(q-1)$.
Denote the class of $z\in k^\times$ modulo $U^2$ by $[z]$.
The following collision criterion will also be used for $q=3$.

\begin{lemma}\label{lem:higher-collision}
Suppose $\gcd(r,q-1)=1$ and
$r\ell_t\equiv2\pmod E$ for some $1\leqslant t<u$.
Let $\kappa\in k^\times\setminus\mathbb{F}_q^\times$, put
$b=\kappa^{-r}$ and $d=\kappa^{-r-D}$, and suppose that
$b,d\notin\mathbb{F}_q$. Fix $a\in L^\times$.
For each $\lambda\in K^\times$ there is a unique $c\in K^\times$
with $c^r=\lambda$. The number of $x\in L^\times$ satisfying
$F(x)=F(\kappa cx)$ is exactly $D$ if
\[
 \left[N_{L/k}\left(\frac{\lambda-b}{\lambda-d}\right)\right]
       =[N_{L/k}(-a)]^{-1},
\]
and zero otherwise. In particular, the displayed condition yields
a collision on $L^\times$.
\end{lemma}

\begin{proof}
The congruence and oddness of $E$ imply $\gcd(r,Q-1)=1$ and
$c^{r+D}=(c^r)^{q^t}$ on $K^\times$. With $c^r=\lambda$,
the collision equation for $x\ne0$ becomes
\[
       x^D=a\frac{\kappa^r\lambda-1}
                       {1-\kappa^{r+D}\lambda^{q^t}} .
\]
Both factors are nonzero, since $K\cap k=\mathbb F_q$.
As $u$ is odd, $2i+t$ ($0\leqslant i<u$) runs through the
residues modulo $u$. Thus the right side has norm
\[
       N_{L/k}(-a)\,\kappa^{-Du}\,
       N_{L/k}\!\left(\frac{\lambda-b}{\lambda-d}\right).
\]
Since $D\mid q^2-1$, the subgroup $(k^\times)^D$ of the cyclic group
$k^\times$ has order $(q+1)/2$, so $(k^\times)^D=U^2$; in particular
$\kappa^{-Du}\in U^2$. The norm $N_{L/k}\colon L^\times\to k^\times$ is
surjective and maps $(L^\times)^D$ into $(k^\times)^D$, so it induces a
surjection
\[
 L^\times/(L^\times)^D\longrightarrow k^\times/(k^\times)^D=k^\times/U^2.
\]
Both groups are cyclic of order $D$, because $D$ divides $q^2-1$ and
$Q^2-1$; hence this map is an isomorphism. Thus an element of $L^\times$
is a $D$-th power if and only if its norm lies in $U^2$. By the norm
computation above, the right side is a $D$-th power exactly under the
displayed class condition, and then the equation has exactly $D$
solutions, which form a coset of $\mu_D\subset L^\times$; otherwise it
has none. Finally
$\kappa\notin K$ makes each pair of inputs distinct.
\end{proof}

\begin{proof}[Proof of Proposition~\ref{prop:higher-odd}]
Suppose $F$ permutes $L$ for $q\geqslant5$. Then $\gcd(r,q-1)=1$,
and Proposition~\ref{prop:higher-exponent}, after inversion if necessary,
gives the required exponent congruence. Take $\kappa$
generating $k^\times$. Since $r$ is odd, $b,d\notin\mathbb F_q$;
$b=d$ would imply $q+1\mid2$, while $b^q=d$ would imply
$q+1\mid r-2$, impossible by parity. Hence the four conjugate
points are distinct. Put $g=[N_{L/k}(-a)]^{-1}$.
If $u=3$ and $g=1$, then $-a$ is a $D$-th power and $F$
has a nonzero root. For $u\geqslant5$, Proposition~\ref{prop:higher-norm}
provides at least $u>2$ points of class $g$, including some
$\lambda\in K^\times$. For $u=3$ and $g\ne1$, we have
$N(b)=N(d)=N(\kappa)^{-r}$ nonsquare and
$d/b=\kappa^{-D}\in U^2$; the same proposition supplies a point
of class $g$ in $K^\times$, since $[W_3(0)]=[W_3(\infty)]=1$.
Lemma~\ref{lem:higher-collision} now gives the contradiction.
\end{proof}

\subsection{The field of order three}\label{sec:higher-ternary}

\begin{proposition}\label{prop:higher-ternary}
If $e\geqslant4$ is even, $r\geqslant1$ and $a\in\mathbb{F}_{3^e}^{\times}$,
then $X^r(X^4+a)$ does not permute $\mathbb{F}_{3^e}$.
\end{proposition}
\begin{proof}
By Lemma~\ref{lem:higher-descent} and Proposition~\ref{prop:higher-quartic},
it suffices to consider $e=2u$, $u\geqslant3$ odd. Let
$Q=3^u$, $K=\mathbb{F}_Q$, $L=\mathbb{F}_{Q^2}$, $k=\mathbb{F}_9$, and
$E=(Q-1)/2$. If a permutation existed, the normalization
in Proposition~\ref{prop:higher-exponent} would give $a=\gamma\in k$
of order four, $r$ odd, and
$r\ell_t\equiv2\pmod E$ for some $1\leqslant t<u$.
Take a generator $\kappa$ of $k^\times$ and put
$b=\kappa^{-r}$, $d=\kappa^{-r-4}=-b$.
The class $[N_{L/k}(-\gamma)]=[(-\gamma)^u]$ in
$k^\times/\{\pm1\}$ is the unique element of order two.
Corollary~\ref{cor:higher-norm-three} supplies
$\lambda\in\mathbf P^1(K)$ with the required norm-ratio
class. The points $0,\infty$ have trivial classes, hence
$\lambda\in K^\times$. Lemma~\ref{lem:higher-collision} yields
a collision, a contradiction.
\end{proof}

\subsection{Proof of the nonexistence theorem}\label{subsec:higher-proof}

\begin{proof}[Proof of Theorem~\ref{thm:higher}]
For $q=3$ apply Proposition~\ref{prop:higher-ternary}. Suppose $q\geqslant5$ and choose the
least even $e\geqslant4$ admitting a permutation, if one exists.
Proposition~\ref{prop:higher-quartic} excludes $e=4$.
If $4\mid e$, Lemma~\ref{lem:higher-descent} gives a smaller even degree
$e/2\geqslant4$. Otherwise $e=2u$ with $u\geqslant3$ odd, contrary to
Proposition~\ref{prop:higher-odd}.
\end{proof}

\section{The complete classification and exponent differences \texorpdfstring{$d(q-1)$}{d(q-1)}}\label{sec:gcd-two}

\begin{proof}[Proof of Theorem~\ref{thm:classification}]
For $a=0$, $F_{r,a}=X^{r+2(q-1)}$; for $e=1$, $F_{r,a}=(1+a)X^r$ as a
function on $\mathbb{F}_q$. This gives Parts~\ref{class:zero}
and~\ref{class:base}. In the situation of Part~\ref{class:odd-index},
$\ell_e$ is odd, so $\gcd(2,\ell_e)=1$, and the assertion follows
directly from \cite[Corollary~6.1]{Fan26}, applied with $d=2$.
Parts~\ref{class:higher-even} and~\ref{class:quadratic} are
Theorems~\ref{thm:higher} and~\ref{thm:quadratic}. The parameter ranges
of the five parts exhaust the family.
\end{proof}

We next reduce the exponent difference of $f_{r,d,a}$ without changing
$a$.

\begin{lemma}\label{lem:general-d-reduction}
Let $q$ be a prime power, $e,d\geqslant1$, and put
$\eta=\gcd(d,\ell_e)$. There is a positive integer $t$ such that
$\gcd(t,q^e-1)=1$ and $dt\equiv\eta\pmod{\ell_e}$.
For all $r\geqslant1$ and $a\in\mathbb{F}_{q^e}$,
\[
 f_{r,d,a}(X^t)=f_{rt,\eta,a}(X)
\]
as functions on $\mathbb{F}_{q^e}$. In particular, these two
polynomials have the same permutation property. When $a\ne0$,
$f_{r,d,a}$ has no nonzero root if and only if
$(-a)^{\ell_e/\eta}\ne1$.
\end{lemma}
\begin{proof}
Put $M=q^e-1$, $H=\ell_e/\eta$, and $\delta=d/\eta$.
As $\gcd(\delta,H)=1$, the Chinese remainder theorem permits
$\delta t\equiv1\pmod H$ and $t\equiv1\pmod\pi$ for every
prime $\pi\mid M$ not dividing $H$. This $t$ is prime to $M$,
and $dt\equiv\eta\pmod{\ell_e}$. Hence
$dt(q-1)\equiv\eta(q-1)\pmod M$, proving the identity
on nonzero elements; both sides vanish at zero. The image
of $x\mapsto x^{d(q-1)}$ has order
$M/\gcd(M,d(q-1))=\ell_e/\eta$, proving the root condition.
\end{proof}

\begin{proof}[Proof of Corollary~\ref{cor:gcd-two}]
Since $\ell_e$ is even, $q$ is odd and $e$ is even.
For $a=0$ the assertion is immediate. For $a\ne0$, apply
Lemma~\ref{lem:general-d-reduction} with $\eta=2$ and
Theorem~\ref{thm:classification} to $F_{rt,a}$. If $e=2$, then
$t\equiv\delta^{-1}\pmod h$, and
\[
 [rt-2]_h=j,\qquad
 \gcd(rt,q-1)=\gcd(r,q-1),\qquad
 \gcd(rt-2,q+1)=\gcd(r-d,q+1).
\]
These identities give precisely the stated conditions.
\end{proof}

For $\eta=2$ the nonzero-root condition is
$(-a)^{\ell_e/2}\ne1$, not $(-a)^{\ell_e}\ne1$.
Together with \cite[Corollary~6.1]{Fan26}, this classifies
all nonzero-coefficient cases with $\eta\leqslant2$ and $e\geqslant2$.
For $a=0$ the criterion is $\gcd(r+d(q-1),q^e-1)=1$;
for $e=1$ it is $a\ne-1$ and $\gcd(r,q-1)=1$.

\providecommand{\bysame}{\leavevmode\hbox to3em{\hrulefill}\thinspace}
\providecommand{\MR}{\relax\ifhmode\unskip\space\fi MR }
\providecommand{\MRhref}[2]{%
  \href{http://www.ams.org/mathscinet-getitem?mr=#1}{#2}
}
\providecommand{\href}[2]{#2}
\Needspace{16\baselineskip}


\begin{thebibliography}{NPYY26}

\bibitem[Fan26]{Fan26}
Xiang Fan, \emph{A complete classification of permutation binomials of the form
  {$X^r(X^{q-1}+a)$} over finite fields}, 2026,
  \href{https://arxiv.org/abs/2609.20354v1}{arXiv:2609.20354v1}.

\bibitem[Hou16]{Hou16}
Xiang-dong Hou, \emph{Permutation polynomials of the form
  {$X^r(a+X^{2(q-1)})$}---a nonexistence result}, 2016,
  \href{https://arxiv.org/abs/1609.03662v1}{arXiv:1609.03662v1}.

\bibitem[HP23]{HPL23}
Xiang-dong Hou and Vincenzo {Pallozzi Lavorante}, \emph{New results on
  permutation binomials of finite fields}, Finite Fields Appl. \textbf{88}
  (2023), 102179,
  \href{https://doi.org/10.1016/j.ffa.2023.102179}{doi:10.1016/j.ffa.2023.102179};
  preprint \href{https://arxiv.org/abs/2111.06533v2}{arXiv:2111.06533v2}.

\bibitem[MRS22]{MRS22}
A.~M. Masuda, I.~Rubio, and J.~Santiago, \emph{Permutation binomials of the
  form {$x^r(x^{q-1}+a)$} over {$\mathbb{F}_{q^e}$}}, Finite Fields Appl.
  \textbf{79} (2022), 102003,
  \href{https://doi.org/10.1016/j.ffa.2022.102003}{doi:10.1016/j.ffa.2022.102003}.

\bibitem[NPYY26]{NPYY26}
Junna Ni, Xuan Pang, Jianhua Yu, and Pingzhi Yuan, \emph{Permutation binomials
  of the form {$X^r(X^{q-1}+a)$} over finite fields}, 2026,
  \href{https://arxiv.org/abs/2609.28595v1}{arXiv:2609.28595v1}.

\end{thebibliography}
\end{document}